\documentclass[draft]{agujournal2019}
\usepackage{url} 
\usepackage{soul}

\usepackage{times}
\usepackage{graphicx}
\usepackage{subcaption}
\usepackage{amsmath}
\usepackage{amsfonts}
\usepackage{amssymb}
\usepackage{amsthm}
\usepackage{verbatim}
\usepackage{latexsym}
\usepackage{float,booktabs}
\usepackage{mathtools}
\usepackage{multicol,multirow}
\usepackage[ruled,vlined]{algorithm2e}
\usepackage{cleveref}
\usepackage{array}

\draftfalse
\journalname{}

\numberwithin{equation}{section}
\allowdisplaybreaks

\theoremstyle{plain}

\theoremstyle{remark}

\theoremstyle{remark}

\begin{document}
 \graphicspath{
               {Figures/}
              }

\title{Data-informed Emulators for Multi-Physics 
Simulations
 }


\authors{Hannah Lu\affil{1}, Dinara Ermakova\affil{2}, \\
Haruko Murakami Wainwright\affil{2,3}, Liange Zheng\affil{3}, Daniel M. Tartakovsky\affil{1}}

\affiliation{1}{Department of Energy Resources Engineering, Stanford University, Stanford, CA 94305, USA}
\affiliation{2}{Department of Nuclear Engineering, University of California Berkeley, Berkeley, CA 94720, USA}
\affiliation{3}{Earth and Environmental Sciences Division, Lawrence Berkeley National Laboratory, Berkeley, CA 94720, USA}

\correspondingauthor{Daniel Tartakovsky}{tartakovsky@stanford.edu}




\begin{keypoints}
\item We developed accurate and efficient surrogates for multi-physics simulations
\item The surrogates' accuracy is greatly enhanced by clustering of training data
\item The surrogates can be used in system-level assessment of a nuclear waste repository
\end{keypoints}

%
%

%
%


\begin{abstract}
Machine learning techniques are powerful tools for construction of emulators for complex systems. We explore different machine learning methods and conceptual methodologies, ranging from functional approximations to dynamical approximations, to build such emulators for coupled thermal, hydrological, mechanical and chemical processes that occur near an engineered barrier system in the nuclear waste repository. Two nonlinear approximators, random forests and neural networks, are deployed to capture the complexity of the physics-based model and to identify its most significant hydrological and geochemical parameters. Our emulators capture the temporal evolution of the Uranium distribution coefficient of the clay buffer, and identify its functional dependence on these key parameters. The emulators' accuracy is further enhanced by assimilating relevant simulated predictors and clustering strategy. The relative performance of random forests and neural networks shows the advantage of ensemble learning in random forests algorithm, especially for highly nonlinear problems with limited data.

\end{abstract}


\section{Introduction}

Subsurface models play a significant role in quantitative analysis and optimization of many environmental phenomena, including nuclear waste disposal~\cite{bea2013identifying}, geological CO$_2$ sequestration~\cite{audigane2007two}, geothermal reservoirs~\cite{xiong2013coupled} and subsurface contamination and remediation~\cite{steefel2015reactive}. A recent focus has been to incorporate  multi-scale and multi-physics models, which typically comprise a large number of coupled (nonlinear) ordinary and partial differential equations. A representative example is  thermal-hydrological-mechanical-chemistry (THMC)  models~\cite{steefel2005reactive,steefel2015reactive,rutqvist2014modeling,zheng2017coupled} used to represent, e.g.,  the changes in flow characteristics due to subsurface evolution caused by thermal and chemical processes. Despite continuing advances in software and hardware development, including high-performance computing, multi-physics simulations with large degrees of freedom remain a demanding and elusive task. That is especially so in sensitivity analysis, uncertainty quantification and inverse modeling where many simulation runs are required. 

Model reduction techniques can significantly reduce the (prohibitively) high computational cost of physics-based simulations, while capturing key features of the underlying dynamics. Such techniques have been used extensively in subsurface applications \cite{schmit1974some,barthelemy1993approximation,simpson2001metamodels,lucia2004reduced,saridakis2008soft,forrester2009recent,razavi2012review} and can be grouped in two general classes. The first is physics-based reduced-order models (ROMs), which seek to map a high-dimensional model onto a meaningful representation of reduced dimensionality; in this context, dimensionality refers to the number of degrees of freedom in a discretized numerical model. 
A prime example of this class is proper orthogonal decomposition (POD)~\cite{kerschen2005method,rowley2005model}, which is grounded in singular value decomposition (SVD). It obtains a ROM by projecting the dynamics of the full model onto the hyperplane using the basis extracted from the SVD analysis. The computational saving stems from replacing the high-dimensional full nonlinear system with its  lower-dimensional counterpart for future prediction. Such ROMs are physics-based in the sense that they inherit the dynamic operator from the projection. 

The second class of ROMs are emulators or surrogates. Instead of reducing a model's dimensionality, these methods aim to reduce its complexity by learning the dynamics of the state variables or quantities of interest directly from the full model's output and/or observational data. 
These data-informed and equation-free ROMs are built by using such machine learning techniques as Gaussian process regression~\cite{rasmussen2003gaussian, pau2013reduced}, dynamic mode decomposition (DMD)~\cite{schmid2010dynamic,kutz2016dynamic}, random forest (RF)~\cite{booker2014comparing,naghibi2016gis}, and neural networks (NN)~\cite{hesthaven2018non,qin2019data}. 

Construction of both types of ROMs for multi-physics (e.g., THMC) problems faces several challenges. First, once discretized in space, complex multi-physics problems result in huge systems of nonlinear ordinary differential equations for which the projection-based techniques become unfeasible. Although POD can be combined with the empirical interpolation method \cite{maday2013generalized, chaturantabut2010nonlinear} in order to handle nonlinearities, it still requires one to solve for a large number of state variables from the projected system. Second, the time evolution of a quantity of interest, expressed as a function of the parameters and simulated predictors, is usually highly nonlinear. It becomes challenging for conventional approximators like Gaussian Process and polynomial regression to capture the dynamics precisely. Third, the computational cost of the full model is so high that only limited amount of high-fidelity data are available for training, which poses a great challenge of overfitting. For these and other reasons,  surrogate models for complex coupled processes  \cite{bianchi2016combining} are scarce. 

Driven by the practical considerations mentioned above, we focus on the data-informed/equa\-tion-free emmulators. Specifically, we investigate the performance of random forests (RF)- and neural networks (NN)-based emulators for complex multi-physics problems. These two surrogates aim to directly predict the dynamics of quantities of interest, without having to deal with the full set of state variables and the governing nonlinear equations.  Among the many machine learning tools, we choose RF and NN, because they are known as robust universal nonlinear approximators that place no formal constrains on data. Conventional regression methods, such as polynomial regression and Gaussian process regression, are expected to fail for complex multi-physics problems because the correlation and/or smoothness conditions they place on the data are seldom satisfied. 

Recent theoretical and computational developments in machine learning (e.g., regularization, cross-validation and bootstrap aggregating) enhance the generalizability of RF and NN, enabling them to handle ``small data''. Since nonlinear systems often exhibit dramatic changes in the relationship between parameters and target variables, we utilize clustering tools to identify the threshold behaviors and design more efficient training strategies. Finally, to boost the computational efficiency of the emulator training, we deploy the machine learning toolboxes \texttt{keras} and \texttt{scikit-learn}.

We use the THC model of an engineered-barrier system at a hypothetical nuclear waste disposal site~\cite{ermakova2020} to illustrate the performance of our RF and NN surrogates. While the model, which consists of a large number of coupled partial-differential equations, describes the spatio-temporal evolution of multiple physico-chemical state variables, a quantity of interest is the Uranium distribution coefficient $K_\text{d}$ for the buffer material. The novelty of our study is two-fold. From the methodological prospective, we improve the prediction accuracy of our emulators by  splitting the training based on the geochemical features; The thresholds of the geochemical features are identified by a cluster analysis on training data. Afterwards, we train each cluster to handle highly nonlinear and non-monotonic functions $K_\text{d} = K_\text{d}(t)$ by taking advantages of the RF and NN approximators. From the applications prospective, although RF and NN have been used before as emulators of relatively simple subsurface models~\cite{zhou-2020-markov,booker2014comparing,naghibi2016gis}, we are not aware of their use for such complex multi-physics phenomena as THC. 

In section~\ref{sec:problem}, we provide a brief description of the THC model~\cite{ermakova2020}
and identify a relevant quantity of interest. In section~\ref{sec:method}, we detail the general methodology for construction of the RF and NN emulators.  The accuracy and robustness of these two surrogates in the THC context  are investigated in section~\ref{sec:numerical}. Main conclusions drawn from this study are summarized in section~\ref{sec:concl}.

\section{Problem Formulation}
\label{sec:problem}

We consider multi-physics simulations of a phenomenon that is described by $N_\text{sv}'$ state variables $\mathbf s(\mathbf x,t) = \{s_1,\dots,s_{N_\text{sv}'} \}$, varying in space $\mathbf x \in \mathcal D$ and time $t \in [0,T]$ throughout the simulation domain $\mathcal D$ during the simulation time interval $[0,T]$. The spatiotemporal evolution of these state variables is described by a system of coupled partial-differential equations
\begin{align}\label{eq:genPDF}
    \frac{\partial s_i}{\partial t} = \mathcal N_i(\mathbf s; \mathbf p), \quad (\mathbf x,t) \in \mathcal D \times (0,T]; \qquad i = 1,\dots,N_\text{sv}',
\end{align}
where $\mathcal N_i$ are (nonlinear) differential operators that contain spatial derivatives, and $\mathbf p = \{p_1,\dots,$ $p_{N_\text{par}} \}$ is a set of $N_\text{par}$ parameters that might vary with in space and time $(\mathbf x,t)$ and be dependent on $\mathbf s$. Problems of this kind have to be solved numerically, which requires a discretization of the spatial domain $\mathcal D$ into $N_\text{el}$ elements (or nodes) and the simulation time horizon $[0,T]$ into $N_\text{st}$ time steps. Consequently, the numerical solution of~\eqref{eq:genPDF} gives a  set of $N_\text{sv} = N_\text{sv}'\times N_\text{el}$ disretized state variables, $\mathbf s_{kn} = \mathbf s(\mathbf x_k,t_n)$ with $k =1,\dots,N_\text{el}$ and $n=1,\dots,N_\text{st}$. 

More often than not, this model output has to be post-processed to compute $N_Q$ quantities of interest (QoIs) $\mathbf Q = \{Q_1,\dots,Q_{N_Q} \}$, such that $Q_i = \mathcal M_i(\mathbf s_{kn})$ for $i=1,\dots,N_Q$, $k =1,\dots,N_\text{el}$, and $n=1,\dots,N_\text{st}$. The maps $\mathcal M_i$ can represent, e.g., numerical approximations of the integrals over $\mathcal D$ or a streamline. In any simulation of practical significance, $N_Q \ll N_\text{sv}$ which makes QoIs much easier to visualize and comprehend than the raw output $\mathbf s_{kn}$. Surrogate modeling aims to derive relationships $\mathbf Q = \mathbf Q(\mathbf p; \mathbf x,t)$ directly, bypassing the need to compute $\mathbf s(\mathbf x,t)$ first.

\begin{table}[!htp]
\begin{center}
\begin{tabular}{ p{0.5cm} p{2cm} p{5cm} p{1.5cm}  p{1cm} p{1cm} }
\hline
$\mathbf p$& Parameter ID & Physical meaning& Reference value&min value&max value \\ 
\hline
$p_1$&ilsoh&
Adsorption surface area on illite (cm$^2/$g) in $\log_{10}$ scale - combined parameter for: 

illite strong site adsorption zone surface area - high sorption affinity; 

illite weak site adsorption zone surface area - lower sorption affinity
&$5$&$3$&$6$\\  
$p_2$&smsoh&
Adsorption surface area on smectite (cm$^2/$g) in $\log_{10}$ scale - combined parameter for:  
smectite strong site adsorption zone surface area - high sorption affinity;
smectite weak site adsorption zone surface area - lower sorption affinity
&$5$&$3$&$6$\\  
$p_3$&pH&initial pore water pH& $7.96$&$9$&$7$\\  
$p_4$&ca$^{2+}$&initial Ca$^{2+}$ concentration in $\log_{10}$ scale in the aqueous phase& $-1.66$&$-3$&$-1$\\
$p_5$&smectite&volume fraction of smectite&$0.92$&$0.3$&$0.95$\\  
$p_6$&illite&volume fraction of illite& $0.0001$&$0.01$&$0.2$\\  
$p_7$&calcite&volume fraction of calcite& $0.01$&$0.01$&$0.03$\\  
\hline
\end{tabular}
\end{center}
\caption{Hydrological and geochemical properties of the engineered-barrier system \protect\cite{ermakova2020}.}
\label{table:params}
\end{table}

To make the exposition concrete, we ground our analysis in the THC model for reactive transport of uranium (U) within  an engi\-neered-barrier system at a hypothetical nuclear waste disposal site~\cite{ermakova2020}. The key component of the system is a clay-based buffer surrounding waste canisters, since clay has high sorption capacity for many radionuclides. The clay properties change over time partly due to thermal and hydrological processes, which may reduce the sorption capacity or degrade the barrier function.  The model~\cite{ermakova2020} 
consists of $N_\text{sv}' = 22$ partial-differential equations for the $N_\text{sv}'$ state variables $\mathbf s(\mathbf x,t)$, representing fluid pressure, saturation, temperature, the concentrations of primary species associated with uranium surface complexation, and the concentrations of around 80 geochemical complexes. In the TOUGHREACT simulations~\cite{xu2014toughreact} of this problem, the one-dimensional domain $\mathcal D$ is discretized into $N_\text{el} = 206$ elements, and the simulation time horizon $T = 10^3$--$10^5$~years into $N_\text{st} = 7444$ time steps. Although the THC model can be developed for a single canister, it is not possible to extend this full model for the entire repository.

\begin{figure}[!htp]
\includegraphics[trim = 0 15 0 0 ,clip,width = 1\textwidth]{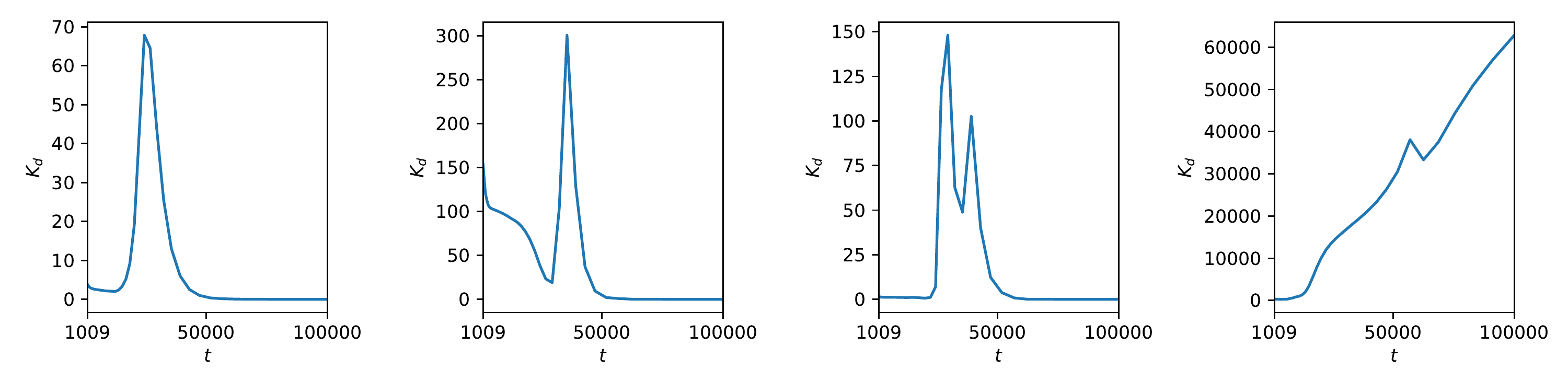}
\caption{Temporal variability of the distribution coefficient $K_\text{d}^{(m)}$, for four combinations of the input parameters $\mathbf p$ indexed by $m$. These realizations of $K_\text{d}$ are evaluated by post-processing the output of the THC model \protect \cite{ermakova2020}. }
\label{fig:example}
\end{figure}

The single QoI from this computation, $Q_1$, is the average distribution coefficient, $K_\text{d}$, of U across the buffer, defined as
\begin{equation}\label{eq:kd}
K_\text{d} = \frac{(\text{total mass of U  sorbed})}{(\text{total mass of U in solution})}.
\end{equation}
This QoI is to be used in a site-scale assessment model. Our goal is to build a surrogate,
\begin{equation}\label{eq:true}
K_{\text{d}} = f(t;\mathbf p), \qquad t \in [0,T], \quad \mathbf p \in \Gamma\subset \mathbb R^{N_\text{par}},
\end{equation}
i.e., to ``learn'' the functional form of $f(\cdot;\cdot)$, from $N_\text{MC}$ solutions of~\eqref{eq:genPDF} obtained for $N_\text{MC}$ different combinations of the parameters $\mathbf p$. The training data are obtained by post-processing the $m$th model run ($m=1,\dots,N_\text{MC}$) to evaluate temporal snapshots (at times $t_1,\dots,t_{N_\text{st}}$) of the corresponding realization of the distribution $K_\text{d}^{(m)}(t)$. Although the THC model~\eqref{eq:genPDF} contains over two hundred  parameters, the model predictions are relatively insensitive to all but seven of them~\cite{ermakova2020}. Therefore, with a slight abuse of notation, we set $N_\text{par} = 7$ in~\eqref{eq:true} while keeping the rest of the model parameters fixed. The seven input parameters used in our examples are collated in Table~\ref{table:params}. The task of learning the function $f(t;\cdot)$ is complicated by the high degree of nonlinearity of $K_\text{d}^{(m)}(t)$ and by the high sensitivity of $K_\text{d}^{(m)}(t)$ to the inputs $\mathbf p$, i.e., by its variability from one value of $m$ to another, (Fig.~\ref{fig:example}).


\begin{figure}[htbp]
\includegraphics[width = 1\textwidth]{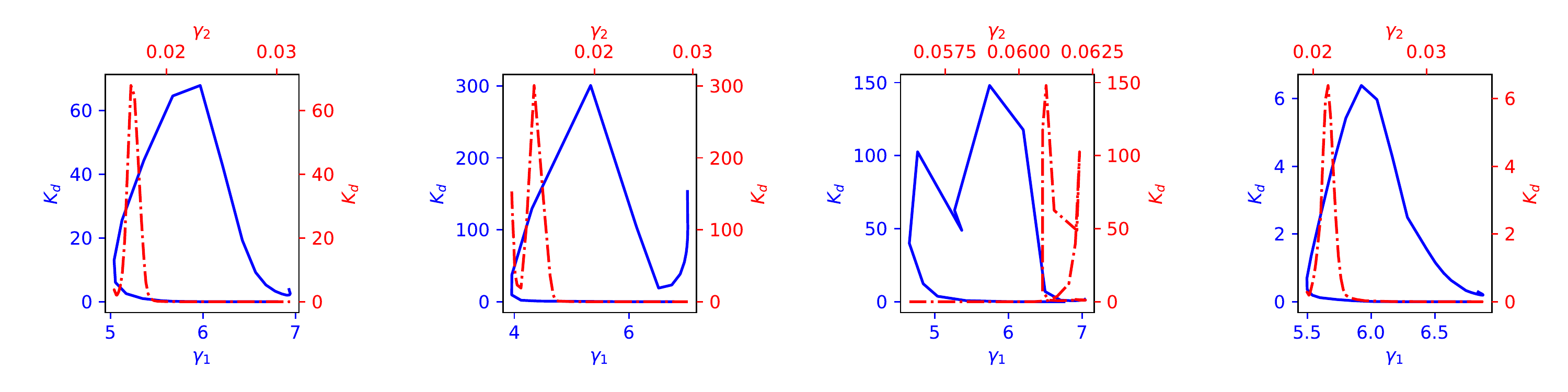}
\caption{Values of $K_\text{d}^{(m)}(t)$ vs. $\gamma_1^{(m)}(t) \equiv \text{pH}(t)$ , and of $K_\text{d}(t)$ vs. $\gamma_2^{(m)}(t) \equiv [\text{Ca}^{2+}](t)$, at the same times $t$, for four combinations of the parameters $\mathbf p$ indexed by $m$. These realizations of $K_\text{d}$ and $\boldsymbol\gamma$ are evaluated by post-processing the output of the THC model.}
\label{fig:gamma_vs_Kd}
\end{figure}

The possible time-dependence of the distribution coefficient $K_\text{d}$ stems from its definition as a map $\mathcal M$ of some of the state variables from the set $\mathbf s(\mathbf x,t)$ or their simulated predictors $\boldsymbol \gamma(\mathbf x,t) = \{\gamma_1, \gamma_2,\dots \}$. In the THC model, these simulated predictors are the pore water composition expressed in term of its pH, $\gamma_1(t) \equiv \text{pH}(t)$ and calcium ion concentration, $\gamma_2(t) \equiv [\text{Ca}^{2+}](t)$, both averaged over the space domain $\mathcal D$. This observation might lead one to attempt to construct an emulator in the form $K_\text{d} = g(\boldsymbol \gamma (t);\mathbf p)$. This formulation is useful when the overall performance model can simulate the regional groundwater chemistry (such as pH and $\text{Ca}^{2+}]$ but not the uranium geochemistry. The observed behavior of $K_\text{d}^{(m)}(t)$ and $\boldsymbol \gamma^{(m)}(t)$ falsifies this hypothesis in all realizations $m$ of the input parameters  (Fig.~\ref{fig:gamma_vs_Kd}): one set of values of $\gamma_1$ and $\gamma_2$ can correspond to two different $K_\text{d}$ values, failing the vertical line tests. This suggests that a surrogate aiming to incorporate the simulated predictors $\boldsymbol \gamma(t)$ must include an explicit dependence on time $t$,
\begin{equation}\label{eq:true2}
K_\text{d} = g(t, \boldsymbol \gamma (t);\mathbf p).
\end{equation}


Another hypothesis is that the present state of $K_\text{d}$ depends not on the present time $t$ but on the whole history of its evolution up to that time. This possible temporal nonlocality of $K_\text{d}$ can be captured by surrogates (see \ref{asec:modified} for details)
\begin{equation}\label{eq:ode}
\frac{\text d K_\text{d} }{ \text dt} = \mathcal F (K_\text{d},t;\mathbf p), \qquad K_\text{d}(t = 0; \mathbf p) = K_\text{d}^0(\mathbf p)
\end{equation}
and 
\begin{equation}\label{eq:ode2}
\frac{\text d K_\text{d} }{ \text dt} = \mathcal G (K_\text{d}, \boldsymbol \gamma(t);\mathbf p), \qquad K_\text{d}(t = 0; \mathbf p) = K_\text{d}^0(\mathbf p).
\end{equation}
If the surrogates~\eqref{eq:true} and~\eqref{eq:true2} are thought of as function approximations of $K_\text{d}$, then~\eqref{eq:ode} and~\eqref{eq:ode2} represent their respective dynamic counterparts. These dynamic approximations aim to learn not only $K_\text{d}$ but also its rate of change at any time $t$. The explicit dependence of $\mathcal F$ on $t$ is, once again, dictated by empirical evidence: for any realization $m$ of the parameters $\mathbf p$, plotting $\text d K_\text{d}^{(m)} / \text dt$ against $K_\text{d}^{(m)}$ at the same times $t$ we found their relation to be multi-valued (similar to Fig.~\ref{fig:gamma_vs_Kd}). A numerical approximation of~\eqref{eq:ode} and~\eqref{eq:ode2} is provided in \ref{asec:modified}. 




The machine-learning techniques used to learn the function $f$ and the functionals $g$, $\mathcal F$, and $\mathcal G$ in~\eqref{eq:true}--\eqref{eq:ode2} are presented below.

\section{Methodology}
\label{sec:method}

The temporal evolution of $K_\text{d}$ is highly nonlinear and can be qualitatively dissimilar for different input parameters $\mathbf p$ (Fig.~\ref{fig:example}). Therefore, we deploy RF (section~\ref{sec:RF}) and NN (section~\ref{sec:NN}) to construct the emulator of $K_\text{d}(t;\mathbf p)$. These machine learning techniques are known to be better nonlinear function approximators than Gaussian-process emulators, polynomial regression,  etc. To substantiate this claim, we demonstrated the poor performance of the Gaussian-process emulator (see Fig.~\ref{fig:GP} in~\ref{asec:GP}). 

\subsection{Data Preprocessing}
\label{sec:prepross}

Our data come from multiple realizations of the THC model~\eqref{eq:genPDF} obtained for different combinations of the seven input parameters $\mathbf p$, i.e., the sample points $\bold p^{(m)}$ are drawn from
\begin{equation}\label{eq:sample}
    p_i^{(m)} \in \left\{p_i^\text{min},\frac{2p_i^\text{min}+p_i^\text{max}}{3},\frac{p_i^\text{min}+2p_i^\text{max}}{3}, p_i^\text{max}\right\}, \quad i \in\{ 1,\cdots, N_\text{par} = 7\}.
\end{equation}
However, some combinations are not physical, leading to false runs in the simulator. For example, output pH values were above 14 or the simulation did not reach the final time step due to lack of convergence. The false runs are dropped from the dataset. Then several parameter samples are drawn by perturbation around the sample pool in~\eqref{eq:sample} to generate more $K_\text{d}^{(m)}$ temporal evolution of desirable shapes. This adjustment aims to enrich the dataset with more balanced temporal variability of $K_\text{d}^{(m)}$.

The preprocessing of these data consist of the following steps.

\subsubsection{Normalization of input and output}

Normalization is a technique often used to prepare the data for machine learning. It is necessary for stable convergence and better accuracy when ``features'', e.g., the input parameters $\mathbf p$ in Table~\ref{table:params}, have vastly different ranges. Therefore, we rescale the parameters $\mathbf p\in \Gamma$ to $\tilde{\mathbf p}\in [-1,1]^{N_\text{par}}$ as the first step of data preparation, i.e., $\tilde p_i = \mathcal R_i(p_i)$ for $i = 1,\dots N_\text{par}$ where $\mathcal R_i$ is the rescaling map. If measurements of the simulated predictors $\boldsymbol\gamma(t)$ are available, we normalize them with their initial values in $\mathbf p$, i.e.,  $\tilde \gamma_1(t) = \mathcal R_3(\gamma_1(t))$ and $\tilde \gamma_2(t) = \mathcal R_4(\gamma_2(t))$. In what follows, we drop the tilde to simplify the notation.

Different values of the input parameters $\mathbf p$ can yield an orders-of-magnitude shift in the range of $K_\text{d}(t)$ (Fig.~\ref{fig:example}). To rescale $K_\text{d}$ and to preserve the positivity, we consider $\ln K_\text{d}$ instead of $K_\text{d}$ in the training and testing.

\subsubsection{Decomposition of parameter space}
\label{sec:decomp}

Despite the nonlinearity of the $K_\text{d}$ time series, one can still discern several distinct parameter regimes, i.e., the subdomains of the ($N_\text{par}=7$)-dimensional parameter space $[-1,1]^{N_\text{par}}$. Low values of the initial pH (parameter $p_3$) create high $K_\text{d}$ values at early reaction times (e.g., Fig.~\ref{fig:example}b). A combination of high illite site density ($p_6$), smectite site density ($p_5$) and initial calcium concentration ($p_4$) with middle range pH ($p_3$) leads to $K_\text{d}$ increasing over the observation time (e.g., Fig.~\ref{fig:example}d). The majority of the $K_\text{d}(t)$ shapes is visually Gaussian (e.g., Figs.~\ref{fig:example}a,c).

\begin{figure}[!htp]
\begin{center}
\includegraphics[width = 0.8\textwidth]{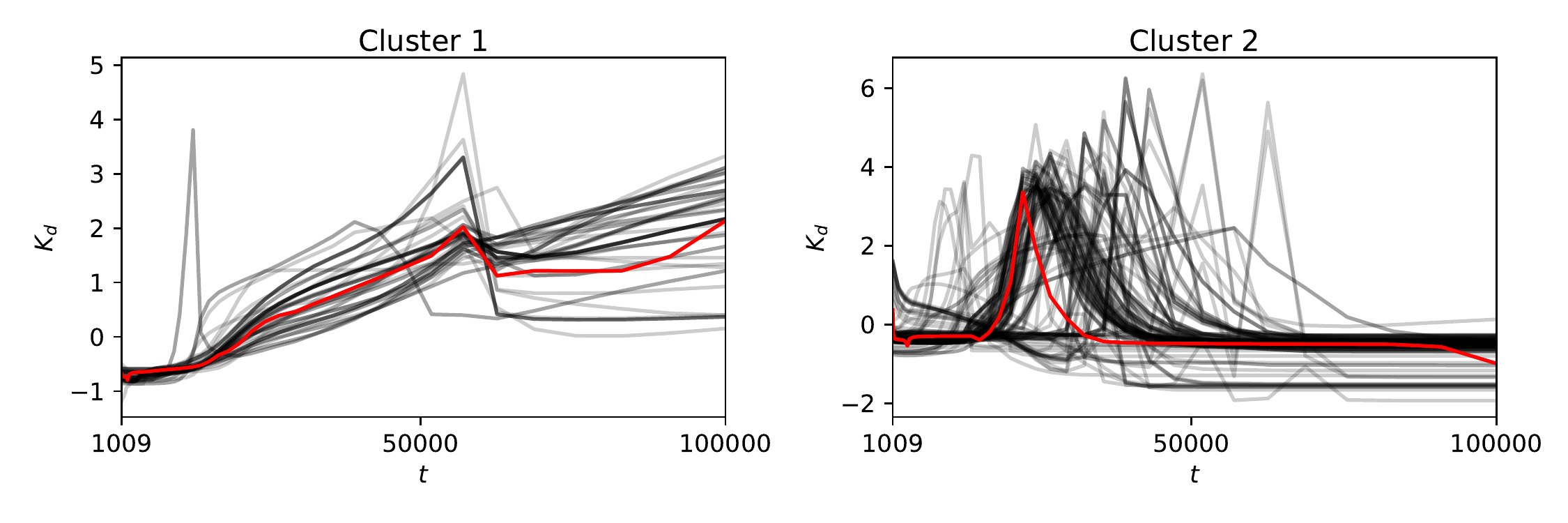}
\end{center}
\caption{Results of classification with $(k=2)$-means clustering with DTW. Each subfigure represents (rescaled) $K_\text{d}^{(m)}$ time series from a given cluster and their centroid (in red).}
\label{fig:cluster}
\end{figure}

Machine-learning tools for classification include  $k$-means, support vector machines, and Gaussian mixtures. Conventional $k$-means techniques perform poorly on time series data, because the Euclidean distance metric is not invariant to time shifts, while most time series data hold such invariants. Therefore, we use $k$-means with dynamic time warping (DTW)~\cite{JMLR:v21:20-091} to deal with time shifts and gather time series of similar shapes. Figure~\ref{fig:cluster} shows the results with $k=2$ clusters. Collating the combinations of the parameters $\mathbf p$ that lead to the $K_\text{d}^{(m)}$ membership in Cluster~1, we 
identify the corresponding parameter subspace 
the region in $[-1,1]^7$ defined by a combination of high $p_2$ and middle range of $p_3$. 

Depending on the input parameter combination, the difference between the lowest and highest output values can be as high as 9 orders of magnitude for Kd due to the wide input parameter ranges. In majority of the cases, depending on the input values of pH ($p_3$) and adsorption surface area on smectite ($p_2$), and a combination of both, Kd may have an increasing shape. A high $p_2$, and 7 $<$ pH $<$ 9 may result in a constant level of Kd or slightly increasing Kd towards the end of the simulation cycle. This results from the increased concentration of bicarbonate ions (HCO3-) at 7 $<$ pH $<$ 9, which leads to the formation of aqueous complexes with U(VI) and high adsorption surface area on smectite reduces the mobility of U(VI) and, as a result, increased Kd - cluster 1. The neutral or acidic (~7.0) or high initial pH and adsorption surface area on smectite lower or equal to 10,000 may lead to an increase of contaminant in aqueous form and a decrease in Kd towards the end of the simulation cycle - cluster 2. 


Based on this clustering observation, we construct the emulators that are trained on each cluster separately. During the test stage, we first determine the test sample's membership in one of the two clusters (based on values of the parameters $p_2$, $p_3$, and $p_4$) and then use the corresponding emulator to predict the temporal evolution of $K_\text{d}$. Our numerical experiments show that the emulators trained on the clustered data have better accuracy than their counterparts trained on the unclustered data.

\subsubsection{Measurement of performance}

The data set consists of time series of $K_\text{d}^{(m)}$ for $N_\text{MC}$ parameter samples $\mathbf p^{(m)}$, with $m=1,\dots,N_\text{MC}$. For each $K_\text{d}^{(m)}$ time series, $M$ snapshots are recorded. The times of these snapshots are logarithmically distributed from $10^3$ years to $10^5$ years in order to better capture the intense reactions at the beginning. To evaluate the performance of the constructed emulators, we reserve $N_\text{test}$ pairs of $\mathbf p^{(m)}$ and the corresponding $K_\text{d}^{(m)}$ time series for validation, while the rest ($N_\text{train} = N_\text{MC} - N_\text{test}$) of the input-output pairs, $\mathbf p^{(m)}$ and $K_\text{d}^{(m)}(t)$, are used for training. Membership in the test set $\mathcal S_{N_\text{test}}$ is determined by randomly drawing $N_\text{test}$ input-output pairs from the total of $N_\text{MC}$ pairs, and the training set is randomly shuffled to reduce the bias caused by sequential ordering of the data. The accuracy of an emulator prediction of the $m$th member from the test set $S_{N_\text{test}}$ is measured by the relative $L_2$ norm,
\begin{equation}\label{eq:rel_error}
\varepsilon_m = \sqrt{
\frac{ \sum_{k = 1}^M\left[\ln K_\text{d}^{(m)}(t_k) - \ln {\tilde K}_\text{d}^{(m)}(t_k) \right]^2 }{  \sum_{k = 1}^M\left[\ln K_\text{d}^{(m)}(t_k) \right]^2 }
}, \qquad \mathbf p^{(m)} \in \mathcal S_{N_\text{test}},
\end{equation}
where $\tilde{K}_\text{d}(t_k;\mathbf p^{(m)})$ is the prediction obtained by the RF or NN emulator. 

\subsection{Random Forest}
\label{sec:RF}

\begin{figure}[!htp]
\begin{center}
\includegraphics[width = 0.9\textwidth]{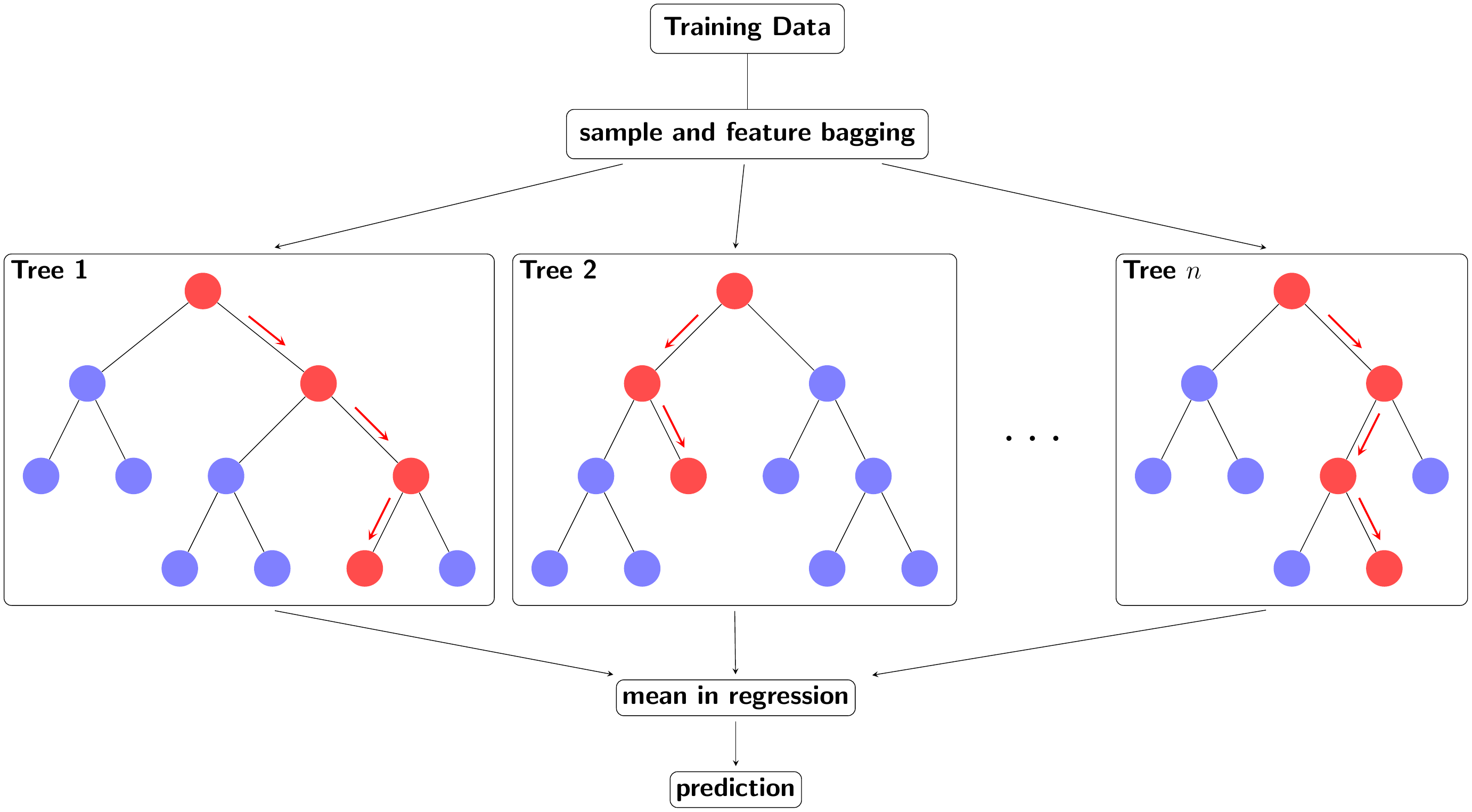}
\end{center}
\caption{Random Forest Learning.}
\label{fig:rf}
\end{figure}

Random Forest (RF) belongs to the group of ensemble learning~\cite{breiman2001random,basu2018iterative}. A group of regression trees~\cite{breiman1984classification} are constructed at training time and the output of RF is the mean prediction of the individual trees, illustrated in Figure~\ref{fig:rf}. Regression trees are known as ``weak learners" in the sense that they have low bias but very high variance, especially for deep trees. In small-data problems, regression trees are seldom accurate due to its high variance~\cite{friedman2001elements}, known as ``overfitting" issue in the machine learning community. To overcome this issue, RF employs multiple regression trees training on different parts of the same training dataset (known as ``sample bagging" strategy~\cite{ho2002data}) so that the average of all trees becomes a ``strong learner", with effectively reduced variance and more accurate learning performance. The input-output data pairs need to be rearranged in the following format:
\begin{equation}\label{eq:rf_input}
\text{input } \mathbf X = [X_1,\cdots,X_s]^\top\in \mathbb R^{s}\to\text{output }Y\in \mathbb R.
\end{equation}
The rank of the importance of each input element $X_i, i = 1,\cdots, s$ is evaluated by \textit{importance scores}~\cite{zhu2015reinforcement} during the fitting process of RF. The details of RF algorithm can be found in the textbook~\cite{friedman2001elements}. We use the RF implementation of the \texttt{RandomForestReg\-ressor} toolbox in \texttt{sklearn} package~\cite{scikit-learn}. Details of implementation are illustrated in section~\ref{sec:numerical}.

\subsection{Neural Networks}
\label{sec:NN}

Artificial neural networks are powerful and robust data-driven modeling tools, especially for nonlinear problems. In conventional notation, the input-output ($\mathbf X$-$\mathbf Y$) maps are approximated by a neural network $\mathcal N_\Theta$:
\begin{equation}
\mathbf Y \approx \mathcal N_\Theta (\mathbf X),\quad \mathbf X\in \mathbb R^s,\quad \mathbf Y\in \mathbb R^r,
\end{equation}
where $\Theta$ is the parameter set including all the parameters in the network. A simplest example is a linear input-output relation $\mathcal N_\Theta = \mathbf W$, where $\mathbf W$ is an $s\times r$ matrix of weights whose numerical values are obtained by minimizing the discrepancy between $\mathbf Y^{(m)}$ and $\mathcal N_\Theta(\mathbf X^{(m)})$, i.e. the following mean squared loss function:
\begin{equation}\label{eq:loss1}
L(\Theta) = \frac{1}{N_\text{train}} \sum_{m = 1}^{N_\text{train}}\|\mathbf Y^{(m)}-\mathcal N_\Theta(\mathbf X^{(m)})\|^2,
\end{equation}
where $\|\cdot \|$ denotes vector $2$-norm hereafter. The performance of this linear regression is likely to be sub-optimal especially for highly nonlinear problems like ours. Thus, one replaces $\mathcal N_\Theta = \mathbf W$ with a nonlinear model $\mathcal N_\Theta = \sigma\circ\mathbf W$, in which the prescribed function $\sigma$ operates on each element of $\mathbf W \mathbf X$. Popular choices for this so-called activation functions include sigmoid, hyperbolic tangent, rectified linear unit (ReLU) and etc. In our numerical tests, we use ReLU as the activation function, i.e., $\sigma(X) = \max (0,X)$. The nonlinear model $\mathcal N_\Theta = \sigma\circ\mathbf W$ constitutes a single fully connected ``layer" in a network.  It has been established that such fully connected NN are universal approximators~\cite{hornik1991approximation,pinkus1999approximation}. 

\begin{figure}[htbp]
\begin{center}
\includegraphics[width = 0.5\textwidth]{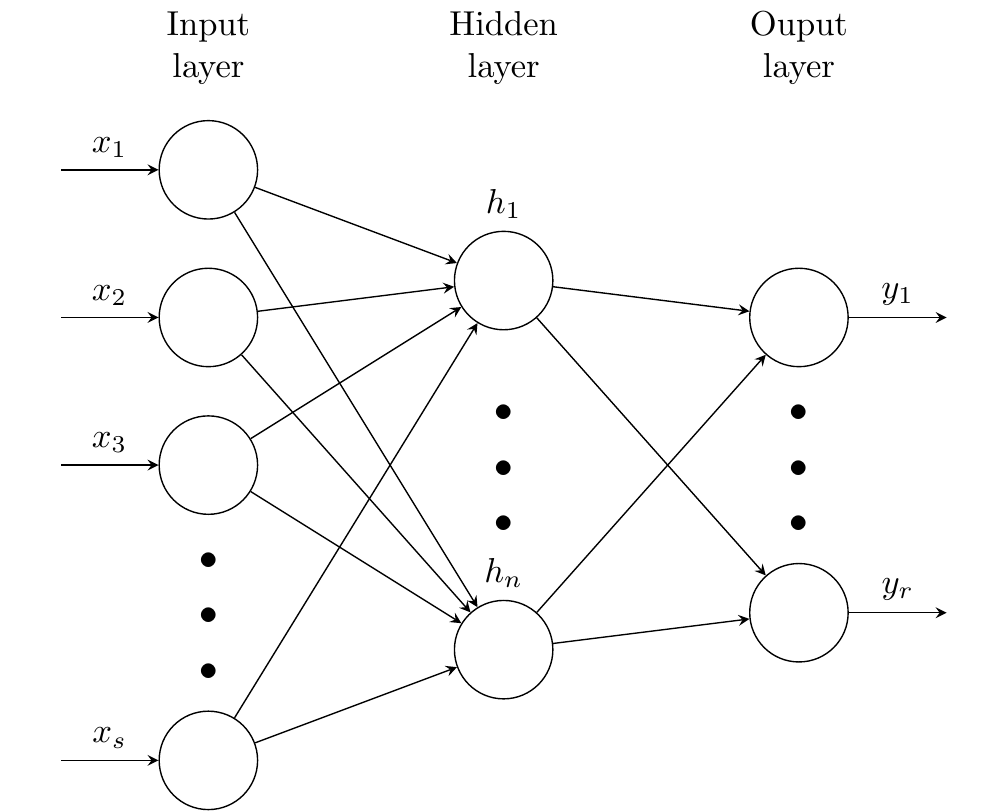}
\end{center}
\caption{An example of $3$-Layer fully connected NN architecture visualization.}
\label{fig:nn}
\end{figure}

A (deep) fully connected NN comprising $L \geq 3$ ``layers" is constructed by a repeated application of the activation function to the input,
\begin{equation}
\mathcal N_\Theta = (\sigma_L\circ \mathbf W_{L-1})\circ \cdots \circ (\sigma_2 \circ \mathbf W_1).
\end{equation}
In general, different activation functions might be used in one network and the last $\sigma_L$ is an identity function, i.e., $\sigma_L(X) = X$. The layers except the input and output layers are called ``hidden layers". The parameter set $\Theta = \{\mathbf W_1,\cdots, \mathbf W_{L-1}\}$ consists of the weights $\mathbf W_l$ connecting the neurons from $l$th to $(l+1)$st layers. The weights $\mathbf W_1$ form a $s\times n_2$ matrix, $\mathbf W_2$ form a $n_3\times n_2$ matrix, $\cdots$, and $W_{L-1}$ is a $n_{L-1}\times r$ matrix, where the integers $n_l, l = 2,\cdots, L-1$ represent the number of neurons in each hidden layer. An example of a $3$-layer fully connected NN architecture is shown in Figure~\ref{fig:nn}. As before, the fitting parameters $\Theta$ are obtained by minimizing~\eqref{eq:loss1}. In practice, $L_2$ norm of the weights is added to the loss function~\eqref{eq:loss1} with small hyperparameter $\lambda$ to avoid overfitting. The learning performance is evaluated by the prediction of $\mathcal N_\Theta(\mathbf X^{(m)})$ compared to $\mathbf Y^{(m)}$ in the test set. NN can be easily implemented using \texttt{TensorFlow Keras} API~\cite{chollet2015keras}. The implementation details are illustrated in section~\ref{sec:numerical}.

\section{Numerical Results and Discussion}
\label{sec:numerical}

Our data set consists of $N_\text{MC} = 172$ input-output pairs, $\{\mathbf p^{(m)}, \boldsymbol{\gamma}^{(m)}(t),K_\text{d}^{(m)}(t) \}_{m=1}^{N_\text{MC}}$. Each output $K_\text{d}^{(m)}(t)$ consists of $M = 50$ snapshots collected at times $t_1,\dots,t_M$. These data set is split into the training and testing data sets consisting of $N_\text{train} = 166$ and $N_\text{test} = 6$ input-output pairs, respectively. The selection of both observation times $\{t_k\}_{k=1}^M$ and membership in the training and testing sets follows the procedures described in section~\ref{sec:prepross}.

The $(k=2)$-means with DTW classifier (section~\ref{sec:decomp}) identifies $N_1 = 123$  and $N_2 = 49$ members in Clusters 1 and 2, respectively (Fig.~\ref{fig:cluster}). Thus, during subgroup learning, we have $N_1^\text{train} = 119$,  $N_1^\text{test} = 4$ and $N_2^\text{train} = 47$,  $N_2^\text{test} = 2$. 

RF is implemented using the machine learning package  \texttt{scikit-learn}. The forest comprises $N_\text{est} = 1000$ regression trees. The maximum depth of each tree is not preset. The nodes are expanded until either all leaves are pure (i.e., cannot be split further) or  all leaves contain less than $1$ sample data. Bootstrap strategy is implemented by drawing the total number of training samples with replacement to fit each tree. In \texttt{RandomForestRegressor} toolbox, the implementation is to set $N_\text{est} = 1000$ with all other default values.

NN is implemented using \texttt{TensorFlow Keras} API. In NN, the neural network consists of $L= 7$ layers with $n_l = n_h$ ($l = 2,\dots,L-1$) neurons in each hidden layers. All the weights are initialized with \texttt{He initialization}~\cite{he2015delving} and all the biases are initialized to be zeros. To avoid overfitting, each layer is penalized by L$_2$ regularization of strength $\lambda$. The training data set is divided into mini-batches of size $10$,  the model is trained for $5000$ epochs, after which the decrease in loss function is saturated. Minimization of the loss function~\eqref{eq:loss1} is done with the \texttt{Adam} algorithm, starting with learning rate $\alpha$. If the monitored validation loss stagnates in a ``patience" $E_p$ epochs, then the learning rate is reduced by a factor of $\beta$ until its preset minimum value. The hyper-parameters $n_h$, $\lambda$, $\alpha$, $\beta$ and the corresponding learning rate schedule need to be fine-tuned and thus are problem-dependent. 

\subsection{Emulators~\eqref{eq:true}: Function Approximation without Observables}
\label{sec:test_true}

To construct the emulator~\eqref{eq:true}, training of the RF~\eqref{eq:rf_input} uses $N = N_\text{train}\times M$ input features $\mathbf X \in \mathbb R^8$ and the same number of output targets $Y \in \mathbb R$,
\begin{align}\label{eq:train_test_data1}
\{\mathbf X^{(i)}\}_{i = 1}^N = \{ t_k, \mathbf p^{(m)} \}_{(m,k)= (1,1)}^{(N_\text{train},M)} \quad\text{and}\quad
\{Y^{(i)}\}_{i = 1}^N = \{\ln K_\text{d}^{(m)}(t_k)\}_{(m,k)= (1,1)}^{(N_\text{train}, M)}.
\end{align}
The testing data-set is arranged identically. The resulting RF-based emulator is denoted by $f_\text{RF}$, such that $\ln K_\text{d}(t) = f_\text{RF}(t; \mathbf p)$. As a curious aside, we found that the importance scores of features $t,p_1,\dots,p_7$, computed by the RF, 
coincide with their rankings obtained via the global sensitivity analysis~\cite{ermakova2020}.


\begin{figure}[!htp]
\begin{center}
\includegraphics[width = 0.9\textwidth]{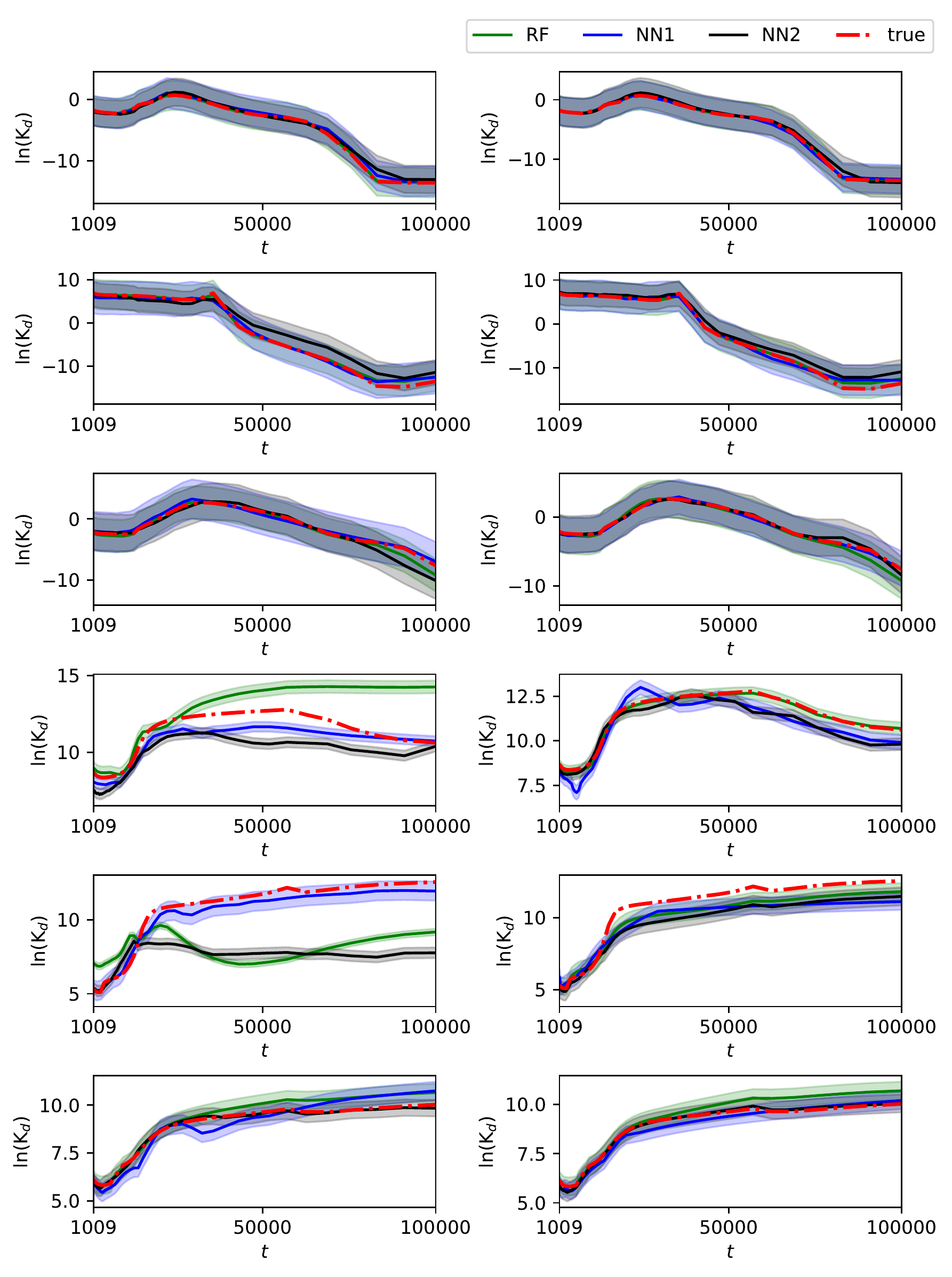}
\end{center}
\caption{RF- and NN-based emulators~\eqref{eq:true} without (left column) and with (right column) data clustering. For the learning with two-class clustering, the top four graphs correspond to Cluster 1 and the remaining bottom two to Cluster 2, with both clusters identified in Fig.~\ref{fig:cluster}. The light green/blue/black region indicates $95\%$ confidence interval.}
\label{fig:results}
\end{figure}

The NN-based emulator~\eqref{eq:true}, trained on the data in~\eqref{eq:train_test_data1}, is denoted by $f_\text{NN1}$ such that $\ln K_\text{d}(t) = f_\text{NN1}(t; \mathbf p)$. As an alternative, we also build an NN-based emulator~\eqref{eq:true} that is trained on $N_\text{train}$ input features $\mathbf X \in \mathbb R^7$ and output targets $\mathbf Y \in \mathbb R^M$,
%
\begin{align}\label{eq:train_test_data2}
\{\mathbf X^{(i)}\}_{i = 1}^{N_\text{train}} = \{ \mathbf p^{(m)} \}_{m=1}^{N_\text{train}}, \qquad
\{\mathbf Y^{(i)}\}_{i = 1}^{N_\text{train}} = \{ \ln K_\text{d}^{(m)}(t_1), \dots, \ln K_\text{d}^{(m)}(t_M) \}_{m=1}^{N_\text{train}}.
\end{align}
We denote this NN-based emulator~\eqref{eq:true} by $f_\text{NN2}$, such that $\{ \ln K_\text{d}(t_1), \dots, \ln K_\text{d}(t_M) \} = f_\text{NN2}(\mathbf p)$. We consider the emulators constructed with and without the clustering of the input data.

The left column of Figure~\ref{fig:results} demonstrates that, when trained on the data without clustering, all the three emulators of type~\eqref{eq:true} yield satisfactory predictions of the distribution coefficient $K_\text{d}(t;\mathbf p)$ in all but a few anomalous cases. This demonstrates the ability of the RF and NN emulators to capture the most common features from the data. The prior clustering of these data enables the  RF and NN emulators to predict the anomalies as well (the right column of Figure~\ref{fig:results}). It leads to significant improvement in the emulators' accuracy.

\begin{table}[htbp]
 \centering
 \begin{tabular}{cclccccccccc}
 \toprule
& \multicolumn{1}{c}{Clustering} 
& \multicolumn{1}{c}{Method} 
& \multicolumn{1}{c}{$i = 1$} 
& \multicolumn{1}{c}{$i = 2$} 
& \multicolumn{1}{c}{$i = 3$} 
& \multicolumn{1}{c}{$i = 4$} 
& \multicolumn{1}{c}{$i = 5$} 
& \multicolumn{1}{c}{$i = 6$}
& \multicolumn{1}{c}{Average} \\
 \toprule
& \parbox[t]{3mm}{\multirow{3}{*}{\rotatebox[origin=c]{90}{No}}}  
& RF & 0.0082 &0.0447& 0.1512 & 0.1107 & 0.2883 & 0.0349 & 0.1063\\
&&NN1&0.0701&0.1015&0.1478 &0.0655&0.0439&0.0517&0.0801\\
 &&NN2&0.1029&0.1717&0.2421&0.1142&0.2745&0.0158&0.1535\\
\hline
 & \parbox[t]{3mm}{\multirow{3}{*}{\rotatebox[origin=c]{90}{Yes}}}   
& RF&0.0055& 0.0433& 0.1595&0.0096&0.0721&0.0365&0.0544\\
 &&NN1&0.0410&0.0700&0.0833&0.0630&0.0900&0.0326&0.0633\\
 &&NN2&0.0706  &0.1212&0.0836&0.0424&0.1019&0.0228&0.0738\\
\bottomrule
 \end{tabular}
 \caption{Relative errors, $\varepsilon_i$, for the $i$th sample and their sample-averages for the RF- and NN-based emulators~\eqref{eq:true} trained on the data without and with clustering.}
 \label{table:error2-2}
 \end{table}
 
In Table~\ref{table:error2-2}, we report the relative errors $\varepsilon_i$ of the three emulators, rendering this assessment more quantitative. Without clustering, NN1 outperforms the other two methods on half of the test cases and on average of all $6$ tests. NN2 performs worst both in most individual cases and in the average evaluation. With clustering, the error in each test cases for every method all drops significantly. On average, RF and NN2 improve their accuracy by nearly $50\%$ respectively. With the help of clustering, RF outperforms the other two methods and provides the most accurate predictions as also shown in Figure~\ref{fig:results}.

\subsection{Emulators~\eqref{eq:true2}: Function Approximation with Observables}
\label{sec:test_true2}

To construct the emulator~\eqref{eq:true2}, training of the RF~\eqref{eq:rf_input} uses $N = N_\text{train}\times M$ input features $\mathbf X \in \mathbb R^{10}$ and the same number of output targets $Y \in \mathbb R$,
\begin{align}\label{eq:train_test_data_m}
\{\mathbf X^{(i)}\}_{i = 1}^N = \{t_k, \boldsymbol\gamma(t_k), \mathbf p^{(m)} \}_{(m,k)= (1,1)}^{(N_\text{train},M)}, \qquad
\{Y^{(i)}\}_{i = 1}^N = \{ \ln K_\text{d}^{(m)}(t_k)\}_{(m,k)= (1,1)}^{(N_\text{train}, M)}. 
\end{align}
The testing data-set is arranged identically. The resulting RF-based emulator is denoted by $f_\text{RF}$, such that $\ln K_\text{d}(t) = g_\text{RF}(t, \boldsymbol\gamma(t); \mathbf p)$.

The NN-based emulator~\eqref{eq:true2}, trained on the data in~\eqref{eq:train_test_data_m}, is denoted by $g_\text{NN1}$ such that $\ln K_{\text{d}} = g_\text{NN1}(t,\boldsymbol\gamma(t);\mathbf p)$. As before, we also construct an alternative NN-based emulator~\eqref{eq:true2} that is trained on $N_\text{train}$ input features $\mathbf X \in \mathbb R^{3M+7}$ and output targets $\mathbf Y \in \mathbb R^M$,
\begin{equation}
\begin{aligned}
&\{\mathbf X^{(i)}\}_{i = 1}^{N_\text{train}} = \{t_1,\dots,t_M, \boldsymbol\gamma(t_1), \dots, \boldsymbol\gamma(t_M), \mathbf p^{(m)} \}_{m=1}^{N_\text{train}},\\
&\{\mathbf Y^{(i)}\}_{i = 1}^{N_\text{train}} = \{ \ln K_\text{d}^{(m)}(t_1), \dots, \ln K_\text{d}^{(m)}(t_M) \}_{m=1}^{N_\text{train}}. 
\end{aligned}
\end{equation}
We denote this NN-based emulator~\eqref{eq:true2} by $g_\text{NN2}$, such that $\{ \ln K_\text{d}(t_1), \dots, \ln K_\text{d}(t_M) \} = g_\text{NN2}(t_1,$ $\dots,t_M, \boldsymbol\gamma(t_1), \dots, \boldsymbol\gamma(t_M); \mathbf p)$. Again, we consider the emulators constructed with and without the clustering of the input data.

\begin{table}[htbp]
 \centering
 \begin{tabular}{cclccccccccc}
 \toprule
& \multicolumn{1}{c}{Clustering} 
& \multicolumn{1}{c}{Method} 
& \multicolumn{1}{c}{$i = 1$} 
& \multicolumn{1}{c}{$i = 2$} 
& \multicolumn{1}{c}{$i = 3$} 
& \multicolumn{1}{c}{$i = 4$} 
& \multicolumn{1}{c}{$i = 5$} 
& \multicolumn{1}{c}{$i = 6$}
& \multicolumn{1}{c}{Average} \\
 \toprule
& \parbox[t]{3mm}{\multirow{3}{*}{\rotatebox[origin=c]{90}{No}}}  
&RF&0.0042&0.0303&0.0660&0.0952&0.0368&0.0026&0.0392\\
 &&NN1&0.0264&0.0643&0.0803&0.1323&0.1416&0.0185&0.0772\\
 &&NN2&0.0497&0.0506&0.1277&0.0898&0.0590&0.0398&0.0694\\
\hline
 & \parbox[t]{3mm}{\multirow{3}{*}{\rotatebox[origin=c]{90}{Yes}}}  
 &RF&0.0037&0.0317&0.0649&0.0335&0.0571&0.0026&0.0322\\
 &&NN1&0.0320&0.0557&0.0593&0.0480&0.0376&0.1187&0.0585\\
 &&NN2&0.0449&0.0527&0.1324&0.0207&0.0605&0.0363&0.0579\\
\bottomrule
 \end{tabular}
 \caption{Relative errors, $\varepsilon_i$, for the $i$th sample and their sample-averages for the RF- and NN-based emulators~\eqref{eq:true2} trained on the data without and with clustering.}
 \label{table:error2-3}
 \end{table}

The performance of these emulators is visually similar to that of their counterparts in  Figure~\ref{fig:results}. The errors $\varepsilon_i$ reported in Table~\ref{table:error2-3} demonstrate that the addition of the simulated predictors, i.e., the use of emulators~\eqref{eq:true2} instead of~\eqref{eq:true}, yields more accurate predictions. The accuracy of each method with simulated predictors and no clustering is already comparable to the accuracy (reported in Table~\ref{table:error2-2}) of each method without simulated predictors but with clustering. This justifies the strong correlations between the simulated predictors and $K_\text{d}$. Similarly as before, clustering helps in improving the accuracy of each method (though less significantly). RF outperforms the other two methods by nearly $50\%$ on average with or without clustering.


\subsection{Emulators~\eqref{eq:ode}: Dynamic Approximation without Observables}
\label{sec:test_ode}

Construction of the RF- and NN-based emulators~\eqref{eq:ode} relies on the same input-output training set. It consists of $N = N_\text{train} \times (M-1)$ input features $\mathbf X \in \mathbb R^{9}$ and output targets $Y \in \mathbb R$,
\begin{equation}\label{eq:train_test_data1m}
\begin{aligned}
&\{\mathbf X^{(i)}\}_{i = 1}^N = \{ \ln K_\text{d}^{(m)}(t_k), t_k, \mathbf p^{(m)} \}_{(m,k)= (1,1)}^{(N_\text{train},M-1)},\\
&\{Y^{(i)}\}_{i = 1}^N = \{ \ln K_\text{d}^{(m)}(t_{k+1}) - \ln K_\text{d}^{(m)}(t_k) \}_{(m,k)= (1,1)}^{(N_\text{train}, M-1)}.
\end{aligned}
\end{equation}
In this case, the NN learning process coincides with \texttt{ResNet}~\cite{he2016deep}. The resulting RF- and NN-based emulators are denoted by $\mathcal F_\alpha$, where $\alpha = \text{RF}$ and $\text{NN}$, respectively (see \ref{asec:modified} for details).

The testing on the data set $\mathcal S_{N_\text{test}}$ is carried out as follows. Given the parameters $\mathbf p^{(m)}$ and the corresponding values $\ln K_\text{d}^{(m)}(t_1)$ from $\mathcal S_{N_\text{test}}$, the values of  $\ln K_\text{d}^{(m)}$ at later times ($t_2,\dots,t_M$) are computed iteratively as
%
\begin{equation}\label{eq:iter1}
\ln K_\text{d}^{(m)}(t_{k+1}) = \ln K_\text{d}^\text{(m)}(t_k) + {\mathcal F}_\alpha( \ln K_\text{d}^\text{(m)}(t_k), t_k; \mathbf p^{(m)}), \qquad
k = 1,\dots,M-1,
\end{equation}
where $\alpha = \text{RF}$ and $\text{NN}$. These predicted values are then compared with the data and the corresponding errors $\varepsilon_m$ are computed.

\begin{table}[htbp]
 \centering
 \begin{tabular}{cclccccccccc}
 \toprule
& \multicolumn{1}{c}{Clustering} 
& \multicolumn{1}{c}{Method} 
& \multicolumn{1}{c}{$i = 1$} 
& \multicolumn{1}{c}{$i = 2$} 
& \multicolumn{1}{c}{$i = 3$} 
& \multicolumn{1}{c}{$i = 4$} 
& \multicolumn{1}{c}{$i = 5$} 
& \multicolumn{1}{c}{$i = 6$}
& \multicolumn{1}{c}{Average} \\
 \toprule
& \parbox[t]{3mm}{\multirow{2}{*}{\rotatebox[origin=c]{90}{No}}}  
&RF&0.0190&0.0921&0.1814&0.0394&0.2390&0.0273&0.0997\\
 &&NN&0.0642&0.1167&0.1109&0.1022&0.1063& 0.0408&0.0902\\
\hline
 & \parbox[t]{3mm}{\multirow{2}{*}{\rotatebox[origin=c]{90}{Yes}}}  
 &RF&0.0087& 0.0686&0.1615&0.0589&0.0448&0.0305&0.0622\\
&&NN&0.0397&0.0792&0.1360&0.0829&0.0429&0.0207&0.0669\\
\bottomrule
 \end{tabular}
 \caption{Relative errors, $\varepsilon_i$, for the $i$th sample and their sample-averages for the RF- and NN-based emulators~\eqref{eq:ode} trained on the data without and with clustering.}
 \label{table:error2-4}
 \end{table}

The predicted distribution coefficients $K_\text{d}^{(m)}(t) \in \mathcal S_{N_\text{test}}$ exhibit the qualitative behavior similar to that in   Figure~\ref{fig:results}. As before, the training with prior clustering improves the the emulators' accuracy. Although the iterative procedure~\eqref{eq:iter1} adds up the error at every time step, one still maintains good accuracy. The errors $\varepsilon_m$ in Table~\ref{table:error2-4} show that the performances of RF and NN are about the same accuracy in the test average error. Comparing Table~\ref{table:error2-2} with Table~\ref{table:error2-4}, we observe that emulator~\eqref{eq:true} and emulator~\eqref{eq:ode} have comparable learning performances.

\subsection{Emulators~\eqref{eq:ode2}: Dynamic Approximation with Observables}
\label{sec:test_ode2}

Construction of the RF- and NN-based emulators~\eqref{eq:ode2} relies on the same input-output training set. It consists of $N = N_\text{train} \times (M-1)$ input features $\mathbf X \in \mathbb R^{12}$ and output targets $Y \in \mathbb R$,
\begin{equation}\label{eq:train_test_data1m}
\begin{aligned}
\{\mathbf X^{(i)}\}_{i = 1}^N =& \{ \ln K_\text{d}^{(m)}(t_k), \boldsymbol\gamma^{(m)}(t_k), \boldsymbol\gamma^{(m)}(t_{k+1}) - \boldsymbol \gamma ^{(m)}(t_k), \mathbf p^{(m)} \}_{(m,k) = (1,1)}^{(N_\text{train},M-1)} \\
\{Y^{(i)}\}_{i = 1}^N =& \{ \ln K_\text{d}^{(m)}(t_{k+1} ) - \ln K_\text{d}^{(m)}(t_k)\}_{(m,k)= (1,1)}^{(N_\text{train}, M-1)}.
\end{aligned}
\end{equation}
The resulting RF- and NN-based emulators are denoted by $\mathcal G_\alpha$, where $\alpha = \text{RF}$ and $\text{NN}$, respectively (see \ref{asec:modified} for details).

The testing on the data set $\mathcal S_{N_\text{test}}$ follows the procedure described in the previous section. Given the parameters $\mathbf p^{(m)}$, the corresponding set of observalbles $\boldsymbol\gamma^{(m)}(t_1),\dots,\boldsymbol\gamma^{(m)}(t_M)$, and the corresponding values $\ln K_\text{d}^{(m)}(t_1)$ from $\mathcal S_{N_\text{test}}$, the values of $\ln K_\text{d}^{(m)}$ at later times ($t_2,\dots,t_M$) are computed iteratively as
\begin{equation}\label{eq:iter2}
\begin{aligned}
\ln K_\text{d}^\text{(m)}(t_{k+1}) = & \ln K_\text{d}^\text{(m)}(t_{k}) + \mathcal G_\alpha( \ln K_\text{d}^\text{(m)}(t_{k}), \boldsymbol \gamma^{(m)}(t_k), \boldsymbol\gamma^{(m)}(t_{k+1}) - \boldsymbol \gamma^{(m)}(t_k); \mathbf p^{(m)}),
\end{aligned}
\end{equation}
for $k=1,\dots,M-1$. These predicted values are then compared with the data and the corresponding errors $\varepsilon_m$ are computed.

\begin{figure}[!htp]
\begin{center}
\includegraphics[width = 0.9\textwidth]{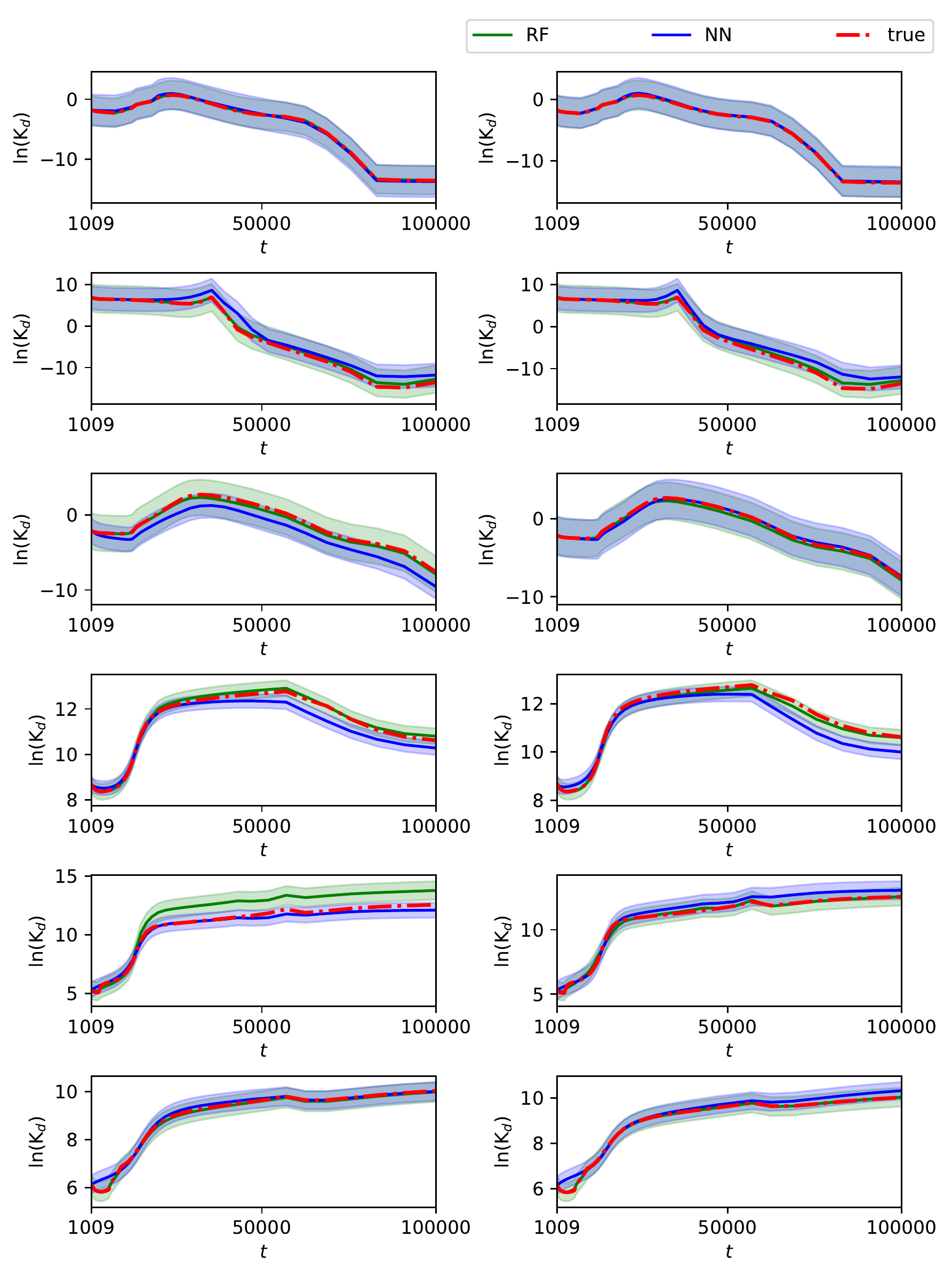}
\end{center}
\caption{RF- and NN-based emulators~\eqref{eq:ode2} without (left column) and with (right column) data clustering. For the learning with two-class clustering, the top four graphs correspond to Cluster 1 and the remaining bottom two to Cluster 2, with both clusters identified in Fig.~\ref{fig:cluster}. The light green/blue/black region indicates $95\%$ confidence interval.}
\label{fig:resnet9}
\end{figure}

The predicted distribution coefficients $K_\text{d}^{(m)}(t) \in \mathcal S_{N_\text{test}}$ are shown in Figure~\ref{fig:resnet9}. The true $K_\text{d}$ lies in the $95\%$ confidence interval of either RF or NN, even with no clustering. The accuracy is further enhanced by the prior clustering, as more clearly shown in the quantitative error evaluations in Table~\ref{table:error2-5}. Although the improvement by using clustering is not significant due to the relatively high accuracy of the emulators without clustering, RF-based emulators~\eqref{eq:ode2} with prior clustering provides the most accurate predictions among all.

\begin{table}[htbp]
 \centering
 \begin{tabular}{cclccccccccc}
 \toprule
& \multicolumn{1}{c}{Clustering} 
& \multicolumn{1}{c}{Method} 
& \multicolumn{1}{c}{$i = 1$} 
& \multicolumn{1}{c}{$i = 2$} 
& \multicolumn{1}{c}{$i = 3$} 
& \multicolumn{1}{c}{$i = 4$} 
& \multicolumn{1}{c}{$i = 5$} 
& \multicolumn{1}{c}{$i = 6$}
& \multicolumn{1}{c}{Average} \\
 \toprule
& \parbox[t]{3mm}{\multirow{2}{*}{\rotatebox[origin=c]{90}{No}}}  
&RF&0.0062 &0.0409&0.0851&0.0079&0.0912&0.0051&0.0394\\
 &&NN&0.0502&0.1504&0.4108&0.0217&0.0278& 0.0330&0.1157\\
\hline
 & \parbox[t]{3mm}{\multirow{2}{*}{\rotatebox[origin=c]{90}{Yes}}}  
&RF&0.0053&0.0513& 0.0867&0.0090&0.0179&0.0032&0.0289\\
 &&NN&0.0270&0.1319&0.0769&0.0287&0.0378&0.0400&0.0570\\
\bottomrule
 \end{tabular}
 \caption{Relative errors, $\varepsilon_i$, for the $i$th sample and their sample-averages for the RF- and NN-based emulators~\eqref{eq:ode2} trained on the data without and with clustering.}
 \label{table:error2-5}
 \end{table}

\subsection{Inter-Method Comparison}

The numerical experiments reported in sections~\ref{sec:test_true}--\ref{sec:test_ode2} demonstrate a comparable performance of the RF- and NN-based emulators, with the former having better accuracy for some test samples and the latter for the others. The data clustering strategy improves the accuracy in all scenarios. 
We summarize the observations and comparisons from the following three perspectives.

\paragraph*{RF- and NN-based emulators.} When no data clustering is performed, the RF-based emulators outperform the NN-based ones in terms of average accuracy. This is due to the power of ensemble learning strategy in the RF algorithm. NNs are designed for big-data problems and, under those condition, perform well as  surrogates \cite{tripathy2018deep,zhu2019physics}. For small-data problems like ours, high variance in a single NN emulator can comprise the prediction accuracy. One can improve the NN emulator by borrowing ideas from RF to construct ensemble neural networks \cite{krogh1995neural,izmailov2018averaging}. This involves training multiple NN models, instead of a single NN model, and combining the predictions from these models. Studies on ensemble learning~\cite{opitz1999popular,polikar2006ensemble,rokach2010ensemble} show that this strategy not only reduces the variance of predictions but also can result in predictions that are better than any single model.

\paragraph*{The NN1 and NN2 designs of NN-based emulators.} The NN1 consider a scalar prediction with time as a input variable, while NN2 consider a time series vector as the target variable. These two input-output designs of NNs  approximate different mappings. In the absence of theoretical explanation of which strategy is to be preferred, we relied on numerical experimentation. We found NN1 to be more stable, i.e., its prediction accuracy is not much affected by the random initialization or train-test splitting. This can be attributed to two factors. First. The input-output design of NN1 enables more data than NN2 ($M$ times more data). Hence, the NN2 prediction may have higher variance and thus be less stable. Second, the NN2 formulation does not consider the correlation between the $t$ and $K_\text{d}(t)$ time series. 
This shortcoming can be overcome by the introduction of input simulated predictors, which account for time correlations, as we have done in sections~\ref{sec:test_true2} and~\ref{sec:test_ode2}. Also, NN1 has better applicability when the emulator is used by another model in which time is usually an input variable.

\paragraph*{Function approximation vs. dynamic formulation.} The RF- and NN-based emulators based on function approximation (sections~\ref{sec:test_true} and~\ref{sec:test_true2}) provide a local representation of the distribution coefficient $K_d(t)$. Their counterparts based on dynamic formulation (sections~\ref{sec:test_ode} and~\ref{sec:test_ode2}) are nonlocal, i.e., account for the memory effect. Since the ``true'' evolution of $K_\text{d}$ is unknown, one cannot determine the preferred formulation without further investigating the dominant dynamics of the upscaled system. Our numerical experiments demonstrate comparable prediction results in both formulations. Other data sets, generated at different scales, may show  different relative performance. Although dynamics formulation would be favorable for real-time estimation, it might not make a significant difference since our predictions are long-term, and direction observations over time is not possible. 

\paragraph*{Impact of simulated predictors.} The incorporation of simulated predictors, e.g., time series of pH and $[\text{Ca}^{2+}]$ used in our examples, significantly improve the learning performance. Since the simulated predictors depend on the geochemical and transport processes, they provide contain information about the full multi-physics system. Therefore, measurements of most correlated simulated predictors play a key role in forecasting upscaled quantities of interest.


\section{Conclusions}
\label{sec:concl}

Development of reduced-order models for reactive transport models remains an open challenge due to nonlinearity and parameter interactions. We constructed RF- and NN-based emulators to represent the buffer-averaged distribution coefficient $K_\text{d}$ as a function of input parameters $\mathbf p$ and time $t$. To the best of our knowledge, ours are the first successful emulators for such systems that are both accurate and computationally efficient. 

We 
explored two formulations of the RF- and NN-based emulators of $K_\text{d}(t;\mathbf p)$: function approximation and dynamic approximation.  We also introduced two strategies to boost the learning performance of all the emulators considered. The first relies on $k$-means clustering with dynamic time warping of the temporal data. 
The second incorporates geochemical simulated predictors, e.g., time series of pore-water pH and calcium concentration, into the learning process. 

The emulators provide orders-of magnitude computational speed-up: the average simulation time for one run of the reactive transport model is about 26 hours, while the training time for a well designed RF- or NN-based emulator is within 10 minutes (on a machine with Intel(R) Core(TM) i7-6700 at 3.40 GHz processor). 

Our use of an emulator can be thought of as numerical upscaling, which is distinct from theoretical upscaling \cite{korneev2016sequential,neuman2009perspective,lichtner2003stochastic}. A good emulator provides not only a meaningful representation of a complex system but also a bridge connecting  high-fidelity (e.g., pore-scale or fine-resolution) models and their low-fidelity (e.g., field-scale) counterparts. In the future, we will implement our $K_\text{d}$ emulators in a performance assessment (PA) model. That effort would consist of the use of the temporally variable $K_\text{d}$ of U(VI) as input parameters in PFLOTRAN simulations of the PA model. In addition, we will test our surrogate models further, improve their accuracy, build in uncertainty quantification, and account for the interactions among parameters.

\section*{Acknowledgements}

Funding for this work was provided by the Spent Fuel and Waste Science and Technology, Office of Nuclear Energy, of the U.S. Department of Energy under Contract Number DE-AC02-05CH11231 with Lawrence Berkeley National Laboratory. Research at Stanford was supported in part by Air Force Office of Scientific Research under award number FA9550-18-1-0474 and the Advanced Research Projects Agency-Energy (ARPA-E), U.S. Department of Energy, under Award Number DE-AR0001202; and a gift from Total. There are no data sharing issues since all of the numerical information is provided in the figures produced by solving the equations in the paper.

\begin{appendix}

\section{Numerical Approximations of $K_\text{d}$ Dynamics}
\label{asec:modified}

Our goal is to provide an accurate approximation to the true solutions of~\eqref{eq:ode} and~\eqref{eq:ode2} at prescribed times $\{t_1,\dots,t_M\}$. Since the analysis of~\eqref{eq:ode} follows directly from that of ~\eqref{eq:ode2}, we show the latter in detail and provide only the final result for the former.

Without loss of generality, we use constant time step
\begin{equation}
\Delta t = t_{k+1}-t_k, \quad \text{for} \quad k = 0,\dots, M-1. 
\end{equation}
For each time interval $[t_k,t_{k+1}]$, with $k = 0,\dots, M-1$, we first seek a first-order local parameterization for the simulated predictors $\boldsymbol \gamma(t)$:
\begin{equation}\label{eq:local_param}
\boldsymbol\gamma(t_k+\tau;\mathbf p)\approx \boldsymbol\gamma(t_k;\mathbf p)+[\boldsymbol\gamma(t_{k+1};\mathbf p)-\boldsymbol\gamma(t_{k};\mathbf p)]\tau, \quad \tau\in[0,\Delta t].
\end{equation}
Then a global parameterization is constructed as
\begin{equation}\label{eq:global_param}
\boldsymbol \gamma(t;\mathbf p)\approx \hat{\boldsymbol \gamma}(t;\mathbf p) \equiv \sum_{k = 0}^{M-1}[\boldsymbol\gamma(t_k;\mathbf p)+(\boldsymbol\gamma(t_{k+1};\mathbf p)-\boldsymbol\gamma(t_{k};\mathbf p))(t-t_k)]\mathbb I_{[t_k,t_{k+1}]}(t),
\end{equation}
where the indicator function $\mathbb I$ is defined by
\begin{equation}
\mathbb I_{[t_k,t_{k+1}]}(t) = \left\{
\begin{aligned}
&1&&\text{if  }t\in[t_k,t_{k+1}],\\
&0&&\text{otherwise.}
\end{aligned}
\right.
\end{equation}
If the simulated predictors $\boldsymbol \gamma(t;\mathbf p)$ are measured either continuously or at more than the two endpoints during all time intervals $[t_k,t_{k+1}]$, then higher-order local parameterizations can be constructed \cite{qin2020data}.

Substituting~\eqref{eq:global_param} into~\eqref{eq:ode2} and using a first-order approximation of the derivative in the latter, we obtain
\begin{equation}\label{eq:ode_soln2_m}
K_\text{d}(t_{k+1};\mathbf p) = K_\text{d}(t_k;\mathbf p) + \Delta t\hat{\mathcal G}( K_\text{d}(t_k), \hat{\boldsymbol \gamma}(t_k),\hat{\boldsymbol \gamma}(t_{k+1})-\hat{\boldsymbol \gamma}(t_k);\mathbf p).
\end{equation}
for $k = 0,\dots, M-1$. Here the target of our approximation $\hat{\mathcal G}: \mathbb R\times \mathbb R^2\times \mathbb R^2\times \mathbb R^7 \to \mathbb R$ is a fully discretized numerical evaluation of $\mathcal G$ with the choice of local parameterization~\eqref{eq:local_param}. 

Similarly, \eqref{eq:ode} is approximated with
\begin{equation}\label{eq:ode_soln_m}
K_\text{d}(t_{k+1};\mathbf p) = K_\text{d}(t_k;\mathbf p) + \hat{\mathcal F}( K_\text{d}(t_k), t_k,\Delta t;\mathbf p),\quad k = 0,\dots, M-1.
\end{equation}
The $\Delta t$-flow maps $\hat {\mathcal F}$ and $\hat {\mathcal G}$ are the target functions for our RF and NN learning methods. The emulators for these flow maps, denoted by $\mathcal G_\alpha$ and  $\mathcal F_\alpha$ (with $\alpha = \text{RF}, \text{NN}$) in~\eqref{eq:iter1} and~\eqref{eq:iter2}, are employed iteratively as surrogates to approximate K$_\text{d}$ values at times $\{t_1,\dots, t_M\}$.

\section{GP Performance}
\label{asec:GP}

\begin{figure}[!htp]
\includegraphics[trim = 0 10 0 0 ,clip,width = 1\textwidth]{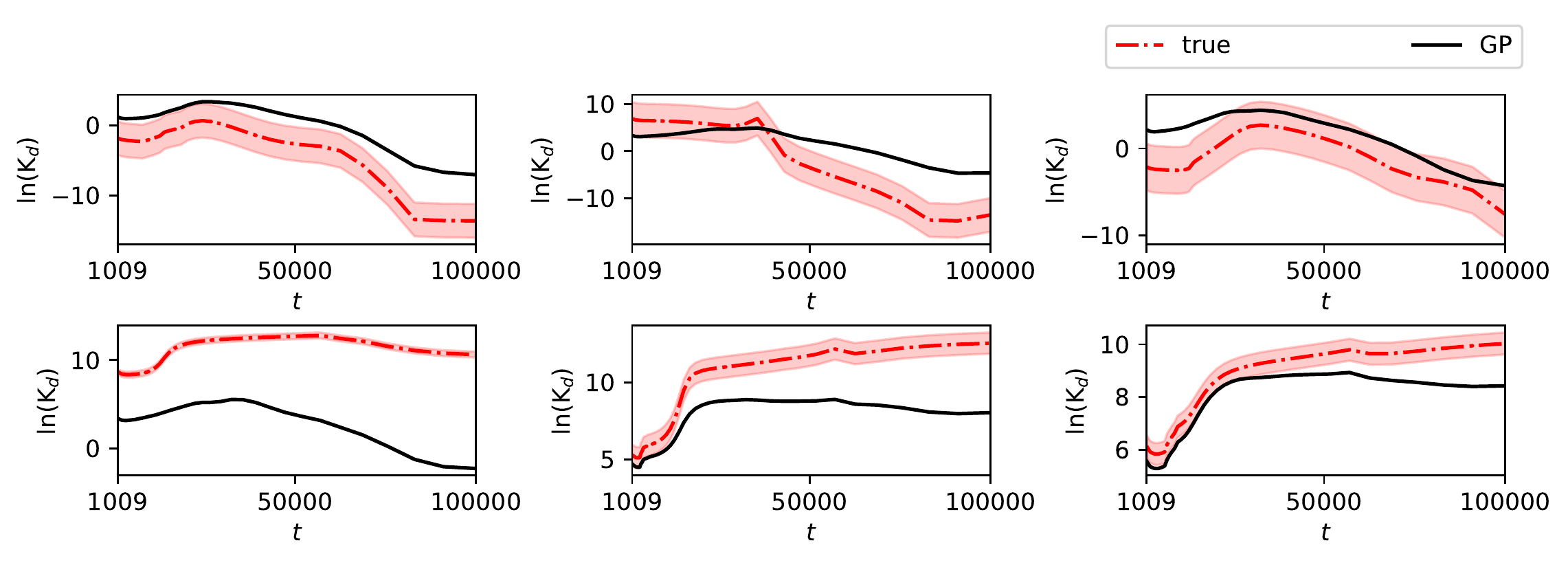}
\caption{GP prediction results: optimal RBF kernel with correlation length $0.522$. The light red region indicates $95\%$ confidence interval.}
\label{fig:GP}
\end{figure}

We use a vanilla GP with common kernels provided in the Python subroutine \texttt{sklearn}. Figure~\ref{fig:GP} demonstrates that it yields the  predictions of $K_\text{d}(t)$ that fall outside the $95\%$ confidence intervals. While more advanced GP variants may improve the learning performance, their exploration lies outside the scope of this work.

\end{appendix}

\bibliography{THC_ROM}

\begin{thebibliography}{}

\bibitem [\protect \citeauthoryear {%
Audigane%
, Gaus%
, Czernichowski-Lauriol%
, Pruess%
\BCBL {}\ \BBA {} Xu%
}{%
Audigane%
\ \protect \BOthers {.}}{%
{\protect \APACyear {2007}}%
}]{%
audigane2007two}
\APACinsertmetastar {%
audigane2007two}%
\begin{APACrefauthors}%
Audigane, P.%
, Gaus, I.%
, Czernichowski-Lauriol, I.%
, Pruess, K.%
\BCBL {}\ \BBA {} Xu, T.%
\end{APACrefauthors}%
\unskip\
\newblock
\APACrefYearMonthDay{2007}{}{}.
\newblock
{\BBOQ}\APACrefatitle {Two-dimensional reactive transport modeling of {CO2}
  injection in a saline aquifer at the {Sleipner} site, {North Sea}}
  {Two-dimensional reactive transport modeling of {CO2} injection in a saline
  aquifer at the {Sleipner} site, {North Sea}}.{\BBCQ}
\newblock
\APACjournalVolNumPages{American Journal of Science}{307}{7}{974--1008}.
\PrintBackRefs{\CurrentBib}

\bibitem [\protect \citeauthoryear {%
Barthe\-lemy%
\ \BBA {} Haftka%
}{%
Barthe\-lemy%
\ \BBA {} Haftka%
}{%
{\protect \APACyear {1993}}%
}]{%
barthelemy1993approximation}
\APACinsertmetastar {%
barthelemy1993approximation}%
\begin{APACrefauthors}%
Barthe\-lemy, J\BHBI F\BPBI M.%
\BCBT {}\ \BBA {} Haftka, R\BPBI T.%
\end{APACrefauthors}%
\unskip\
\newblock
\APACrefYearMonthDay{1993}{}{}.
\newblock
{\BBOQ}\APACrefatitle {Approximation concepts for optimum structural design---a
  review} {Approximation concepts for optimum structural design---a
  review}.{\BBCQ}
\newblock
\APACjournalVolNumPages{Structural Optimization}{5}{3}{129--144}.
\PrintBackRefs{\CurrentBib}

\bibitem [\protect \citeauthoryear {%
Basu%
, Kumbier%
, Brown%
\BCBL {}\ \BBA {} Yu%
}{%
Basu%
\ \protect \BOthers {.}}{%
{\protect \APACyear {2018}}%
}]{%
basu2018iterative}
\APACinsertmetastar {%
basu2018iterative}%
\begin{APACrefauthors}%
Basu, S.%
, Kumbier, K.%
, Brown, J\BPBI B.%
\BCBL {}\ \BBA {} Yu, B.%
\end{APACrefauthors}%
\unskip\
\newblock
\APACrefYearMonthDay{2018}{}{}.
\newblock
{\BBOQ}\APACrefatitle {Iterative random forests to discover predictive and
  stable high-order interactions} {Iterative random forests to discover
  predictive and stable high-order interactions}.{\BBCQ}
\newblock
\APACjournalVolNumPages{Proceedings of the National Academy of
  Sciences}{115}{8}{1943--1948}.
\PrintBackRefs{\CurrentBib}

\bibitem [\protect \citeauthoryear {%
Bea%
\ \protect \BOthers {.}}{%
Bea%
\ \protect \BOthers {.}}{%
{\protect \APACyear {2013}}%
}]{%
bea2013identifying}
\APACinsertmetastar {%
bea2013identifying}%
\begin{APACrefauthors}%
Bea, S\BPBI A.%
, Wainwright, H.%
, Spycher, N.%
, Faybishenko, B.%
, Hubbard, S\BPBI S.%
\BCBL {}\ \BBA {} Denham, M\BPBI E.%
\end{APACrefauthors}%
\unskip\
\newblock
\APACrefYearMonthDay{2013}{}{}.
\newblock
{\BBOQ}\APACrefatitle {Identifying key controls on the behavior of an {acidic-U
  (VI) plume} in the {Savannah River Site} using reactive transport modeling}
  {Identifying key controls on the behavior of an {acidic-U (VI) plume} in the
  {Savannah River Site} using reactive transport modeling}.{\BBCQ}
\newblock
\APACjournalVolNumPages{Journal of Contaminant Hydrology}{151}{}{34--54}.
\PrintBackRefs{\CurrentBib}

\bibitem [\protect \citeauthoryear {%
Bianchi%
, Zheng%
\BCBL {}\ \BBA {} Birkholzer%
}{%
Bianchi%
\ \protect \BOthers {.}}{%
{\protect \APACyear {2016}}%
}]{%
bianchi2016combining}
\APACinsertmetastar {%
bianchi2016combining}%
\begin{APACrefauthors}%
Bianchi, M.%
, Zheng, L.%
\BCBL {}\ \BBA {} Birkholzer, J\BPBI T.%
\end{APACrefauthors}%
\unskip\
\newblock
\APACrefYearMonthDay{2016}{}{}.
\newblock
{\BBOQ}\APACrefatitle {Combining multiple lower-fidelity models for emulating
  complex model responses for {CCS} environmental risk assessment} {Combining
  multiple lower-fidelity models for emulating complex model responses for
  {CCS} environmental risk assessment}.{\BBCQ}
\newblock
\APACjournalVolNumPages{International Journal of Greenhouse Gas
  Control}{46}{}{248--258}.
\PrintBackRefs{\CurrentBib}

\bibitem [\protect \citeauthoryear {%
Booker%
\ \BBA {} Woods%
}{%
Booker%
\ \BBA {} Woods%
}{%
{\protect \APACyear {2014}}%
}]{%
booker2014comparing}
\APACinsertmetastar {%
booker2014comparing}%
\begin{APACrefauthors}%
Booker, D.%
\BCBT {}\ \BBA {} Woods, R.%
\end{APACrefauthors}%
\unskip\
\newblock
\APACrefYearMonthDay{2014}{}{}.
\newblock
{\BBOQ}\APACrefatitle {Comparing and combining physically-based and
  empirically-based approaches for estimating the hydrology of ungauged
  catchments} {Comparing and combining physically-based and empirically-based
  approaches for estimating the hydrology of ungauged catchments}.{\BBCQ}
\newblock
\APACjournalVolNumPages{Journal of Hydrology}{508}{}{227--239}.
\PrintBackRefs{\CurrentBib}

\bibitem [\protect \citeauthoryear {%
Breiman%
}{%
Breiman%
}{%
{\protect \APACyear {2001}}%
}]{%
breiman2001random}
\APACinsertmetastar {%
breiman2001random}%
\begin{APACrefauthors}%
Breiman, L.%
\end{APACrefauthors}%
\unskip\
\newblock
\APACrefYearMonthDay{2001}{}{}.
\newblock
{\BBOQ}\APACrefatitle {Random forests} {Random forests}.{\BBCQ}
\newblock
\APACjournalVolNumPages{Machine Learning}{45}{1}{5--32}.
\PrintBackRefs{\CurrentBib}

\bibitem [\protect \citeauthoryear {%
Breiman%
, Friedman%
, Stone%
\BCBL {}\ \BBA {} Olshen%
}{%
Breiman%
\ \protect \BOthers {.}}{%
{\protect \APACyear {1984}}%
}]{%
breiman1984classification}
\APACinsertmetastar {%
breiman1984classification}%
\begin{APACrefauthors}%
Breiman, L.%
, Friedman, J.%
, Stone, C\BPBI J.%
\BCBL {}\ \BBA {} Olshen, R\BPBI A.%
\end{APACrefauthors}%
\unskip\
\newblock
\APACrefYear{1984}.
\newblock
\APACrefbtitle {Classification and Regression Trees} {Classification and
  regression trees}.
\newblock
\APACaddressPublisher{}{CRC press}.
\PrintBackRefs{\CurrentBib}

\bibitem [\protect \citeauthoryear {%
Chaturantabut%
\ \BBA {} Sorensen%
}{%
Chaturantabut%
\ \BBA {} Sorensen%
}{%
{\protect \APACyear {2010}}%
}]{%
chaturantabut2010nonlinear}
\APACinsertmetastar {%
chaturantabut2010nonlinear}%
\begin{APACrefauthors}%
Chaturantabut, S.%
\BCBT {}\ \BBA {} Sorensen, D\BPBI C.%
\end{APACrefauthors}%
\unskip\
\newblock
\APACrefYearMonthDay{2010}{}{}.
\newblock
{\BBOQ}\APACrefatitle {Nonlinear model reduction via discrete empirical
  interpolation} {Nonlinear model reduction via discrete empirical
  interpolation}.{\BBCQ}
\newblock
\APACjournalVolNumPages{SIAM Journal on Scientific
  Computing}{32}{5}{2737--2764}.
\PrintBackRefs{\CurrentBib}

\bibitem [\protect \citeauthoryear {%
Chollet%
\ \protect \BOthers {.}}{%
Chollet%
\ \protect \BOthers {.}}{%
{\protect \APACyear {2015}}%
}]{%
chollet2015keras}
\APACinsertmetastar {%
chollet2015keras}%
\begin{APACrefauthors}%
Chollet, F.%
\BCBT {}\ \BOthersPeriod {.}
\end{APACrefauthors}%
\unskip\
\newblock
\APACrefYearMonthDay{2015}{}{}.
\newblock
\APACrefbtitle {Keras.} {Keras.}
\newblock
\APAChowpublished {\url{https://keras.io}}.
\PrintBackRefs{\CurrentBib}

\bibitem [\protect \citeauthoryear {%
Ermakova%
, Wainwright%
, Zheng%
, Shirley%
\BCBL {}\ \BBA {} Lu%
}{%
Ermakova%
\ \protect \BOthers {.}}{%
{\protect \APACyear {2020}}%
}]{%
ermakova2020}
\APACinsertmetastar {%
ermakova2020}%
\begin{APACrefauthors}%
Ermakova, D.%
, Wainwright, H.%
, Zheng, L.%
, Shirley, I.%
\BCBL {}\ \BBA {} Lu, H.%
\end{APACrefauthors}%
\unskip\
\newblock
\APACrefYearMonthDay{2020}{}{}.
\newblock
{\BBOQ}\APACrefatitle {Global Sensitivity Analysis for {U(VI)} Transport for
  integrating coupled {THC} Models into {PA} model} {Global sensitivity
  analysis for {U(VI)} transport for integrating coupled {THC} models into {PA}
  model}.{\BBCQ}
\newblock

\PrintBackRefs{\CurrentBib}

\bibitem [\protect \citeauthoryear {%
Forrester%
\ \BBA {} Keane%
}{%
Forrester%
\ \BBA {} Keane%
}{%
{\protect \APACyear {2009}}%
}]{%
forrester2009recent}
\APACinsertmetastar {%
forrester2009recent}%
\begin{APACrefauthors}%
Forrester, A\BPBI I\BPBI J.%
\BCBT {}\ \BBA {} Keane, A\BPBI J.%
\end{APACrefauthors}%
\unskip\
\newblock
\APACrefYearMonthDay{2009}{}{}.
\newblock
{\BBOQ}\APACrefatitle {Recent advances in surrogate-based optimization} {Recent
  advances in surrogate-based optimization}.{\BBCQ}
\newblock
\APACjournalVolNumPages{Progress in Aerospace Sciences}{45}{1-3}{50--79}.
\PrintBackRefs{\CurrentBib}

\bibitem [\protect \citeauthoryear {%
Friedman%
, Hastie%
\BCBL {}\ \BBA {} Tibshirani%
}{%
Friedman%
\ \protect \BOthers {.}}{%
{\protect \APACyear {2001}}%
}]{%
friedman2001elements}
\APACinsertmetastar {%
friedman2001elements}%
\begin{APACrefauthors}%
Friedman, J.%
, Hastie, T.%
\BCBL {}\ \BBA {} Tibshirani, R.%
\end{APACrefauthors}%
\unskip\
\newblock
\APACrefYear{2001}.
\newblock
\APACrefbtitle {The elements of statistical learning} {The elements of
  statistical learning}\ (\BVOL~1)\ (\BNUM~10).
\newblock
\APACaddressPublisher{}{Springer series in statistics New York}.
\PrintBackRefs{\CurrentBib}

\bibitem [\protect \citeauthoryear {%
He%
, Zhang%
, Ren%
\BCBL {}\ \BBA {} Sun%
}{%
He%
\ \protect \BOthers {.}}{%
{\protect \APACyear {2015}}%
}]{%
he2015delving}
\APACinsertmetastar {%
he2015delving}%
\begin{APACrefauthors}%
He, K.%
, Zhang, X.%
, Ren, S.%
\BCBL {}\ \BBA {} Sun, J.%
\end{APACrefauthors}%
\unskip\
\newblock
\APACrefYearMonthDay{2015}{}{}.
\newblock
{\BBOQ}\APACrefatitle {Delving deep into rectifiers: {Surpassing} human-level
  performance on imagenet classification} {Delving deep into rectifiers:
  {Surpassing} human-level performance on imagenet classification}.{\BBCQ}
\newblock
\BIn{} \APACrefbtitle {Proceedings of the {IEEE} International Conference on
  Computer Vision} {Proceedings of the {IEEE} international conference on
  computer vision}\ (\BPGS\ 1026--1034).
\PrintBackRefs{\CurrentBib}

\bibitem [\protect \citeauthoryear {%
He%
, Zhang%
, Ren%
\BCBL {}\ \BBA {} Sun%
}{%
He%
\ \protect \BOthers {.}}{%
{\protect \APACyear {2016}}%
}]{%
he2016deep}
\APACinsertmetastar {%
he2016deep}%
\begin{APACrefauthors}%
He, K.%
, Zhang, X.%
, Ren, S.%
\BCBL {}\ \BBA {} Sun, J.%
\end{APACrefauthors}%
\unskip\
\newblock
\APACrefYearMonthDay{2016}{}{}.
\newblock
{\BBOQ}\APACrefatitle {Deep residual learning for image recognition} {Deep
  residual learning for image recognition}.{\BBCQ}
\newblock
\BIn{} \APACrefbtitle {Proceedings of the {IEEE} Conference on Computer Vision
  and Pattern Recognition} {Proceedings of the {IEEE} conference on computer
  vision and pattern recognition}\ (\BPGS\ 770--778).
\PrintBackRefs{\CurrentBib}

\bibitem [\protect \citeauthoryear {%
Hesthaven%
\ \BBA {} Ubbiali%
}{%
Hesthaven%
\ \BBA {} Ubbiali%
}{%
{\protect \APACyear {2018}}%
}]{%
hesthaven2018non}
\APACinsertmetastar {%
hesthaven2018non}%
\begin{APACrefauthors}%
Hesthaven, J\BPBI S.%
\BCBT {}\ \BBA {} Ubbiali, S.%
\end{APACrefauthors}%
\unskip\
\newblock
\APACrefYearMonthDay{2018}{}{}.
\newblock
{\BBOQ}\APACrefatitle {Non-intrusive reduced order modeling of nonlinear
  problems using neural networks} {Non-intrusive reduced order modeling of
  nonlinear problems using neural networks}.{\BBCQ}
\newblock
\APACjournalVolNumPages{Journal of Computational Physics}{363}{}{55--78}.
\PrintBackRefs{\CurrentBib}

\bibitem [\protect \citeauthoryear {%
Ho%
}{%
Ho%
}{%
{\protect \APACyear {2002}}%
}]{%
ho2002data}
\APACinsertmetastar {%
ho2002data}%
\begin{APACrefauthors}%
Ho, T\BPBI K.%
\end{APACrefauthors}%
\unskip\
\newblock
\APACrefYearMonthDay{2002}{}{}.
\newblock
{\BBOQ}\APACrefatitle {A data complexity analysis of comparative advantages of
  decision forest constructors} {A data complexity analysis of comparative
  advantages of decision forest constructors}.{\BBCQ}
\newblock
\APACjournalVolNumPages{Pattern Analysis \& Applications}{5}{2}{102--112}.
\PrintBackRefs{\CurrentBib}

\bibitem [\protect \citeauthoryear {%
Hornik%
}{%
Hornik%
}{%
{\protect \APACyear {1991}}%
}]{%
hornik1991approximation}
\APACinsertmetastar {%
hornik1991approximation}%
\begin{APACrefauthors}%
Hornik, K.%
\end{APACrefauthors}%
\unskip\
\newblock
\APACrefYearMonthDay{1991}{}{}.
\newblock
{\BBOQ}\APACrefatitle {Approximation capabilities of multilayer feedforward
  networks} {Approximation capabilities of multilayer feedforward
  networks}.{\BBCQ}
\newblock
\APACjournalVolNumPages{Neural Networks}{4}{2}{251--257}.
\PrintBackRefs{\CurrentBib}

\bibitem [\protect \citeauthoryear {%
Izmailov%
, Podoprikhin%
, Garipov%
, Vetrov%
\BCBL {}\ \BBA {} Wilson%
}{%
Izmailov%
\ \protect \BOthers {.}}{%
{\protect \APACyear {2018}}%
}]{%
izmailov2018averaging}
\APACinsertmetastar {%
izmailov2018averaging}%
\begin{APACrefauthors}%
Izmailov, P.%
, Podoprikhin, D.%
, Garipov, T.%
, Vetrov, D.%
\BCBL {}\ \BBA {} Wilson, A\BPBI G.%
\end{APACrefauthors}%
\unskip\
\newblock
\APACrefYearMonthDay{2018}{}{}.
\newblock
{\BBOQ}\APACrefatitle {Averaging weights leads to wider optima and better
  generalization} {Averaging weights leads to wider optima and better
  generalization}.{\BBCQ}
\newblock
\APACjournalVolNumPages{arXiv preprint arXiv:1803.05407}{}{}{}.
\PrintBackRefs{\CurrentBib}

\bibitem [\protect \citeauthoryear {%
Kerschen%
, Golinval%
, Vakakis%
\BCBL {}\ \BBA {} Bergman%
}{%
Kerschen%
\ \protect \BOthers {.}}{%
{\protect \APACyear {2005}}%
}]{%
kerschen2005method}
\APACinsertmetastar {%
kerschen2005method}%
\begin{APACrefauthors}%
Kerschen, G.%
, Golinval, J\BHBI c.%
, Vakakis, A\BPBI F.%
\BCBL {}\ \BBA {} Bergman, L\BPBI A.%
\end{APACrefauthors}%
\unskip\
\newblock
\APACrefYearMonthDay{2005}{}{}.
\newblock
{\BBOQ}\APACrefatitle {The method of proper orthogonal decomposition for
  dynamical characterization and order reduction of mechanical systems: an
  overview} {The method of proper orthogonal decomposition for dynamical
  characterization and order reduction of mechanical systems: an
  overview}.{\BBCQ}
\newblock
\APACjournalVolNumPages{Nonlinear Dynamics}{41}{1-3}{147--169}.
\PrintBackRefs{\CurrentBib}

\bibitem [\protect \citeauthoryear {%
Korneev%
\ \BBA {} Battiato%
}{%
Korneev%
\ \BBA {} Battiato%
}{%
{\protect \APACyear {2016}}%
}]{%
korneev2016sequential}
\APACinsertmetastar {%
korneev2016sequential}%
\begin{APACrefauthors}%
Korneev, S.%
\BCBT {}\ \BBA {} Battiato, I.%
\end{APACrefauthors}%
\unskip\
\newblock
\APACrefYearMonthDay{2016}{}{}.
\newblock
{\BBOQ}\APACrefatitle {Sequential homogenization of reactive transport in
  polydisperse porous media} {Sequential homogenization of reactive transport
  in polydisperse porous media}.{\BBCQ}
\newblock
\APACjournalVolNumPages{Multiscale Modeling \& Simulation}{14}{4}{1301--1318}.
\PrintBackRefs{\CurrentBib}

\bibitem [\protect \citeauthoryear {%
Krogh%
\ \BBA {} Vedelsby%
}{%
Krogh%
\ \BBA {} Vedelsby%
}{%
{\protect \APACyear {1995}}%
}]{%
krogh1995neural}
\APACinsertmetastar {%
krogh1995neural}%
\begin{APACrefauthors}%
Krogh, A.%
\BCBT {}\ \BBA {} Vedelsby, J.%
\end{APACrefauthors}%
\unskip\
\newblock
\APACrefYearMonthDay{1995}{}{}.
\newblock
{\BBOQ}\APACrefatitle {Neural network ensembles, cross validation, and active
  learning} {Neural network ensembles, cross validation, and active
  learning}.{\BBCQ}
\newblock
\BIn{} \APACrefbtitle {Advances in Neural Information Processing Systems}
  {Advances in neural information processing systems}\ (\BPGS\ 231--238).
\PrintBackRefs{\CurrentBib}

\bibitem [\protect \citeauthoryear {%
Kutz%
, Brunton%
, Brunton%
\BCBL {}\ \BBA {} Proctor%
}{%
Kutz%
\ \protect \BOthers {.}}{%
{\protect \APACyear {2016}}%
}]{%
kutz2016dynamic}
\APACinsertmetastar {%
kutz2016dynamic}%
\begin{APACrefauthors}%
Kutz, J\BPBI N.%
, Brunton, S\BPBI L.%
, Brunton, B\BPBI W.%
\BCBL {}\ \BBA {} Proctor, J\BPBI L.%
\end{APACrefauthors}%
\unskip\
\newblock
\APACrefYear{2016}.
\newblock
\APACrefbtitle {Dynamic mode decomposition: data-driven modeling of complex
  systems} {Dynamic mode decomposition: data-driven modeling of complex
  systems}.
\newblock
\APACaddressPublisher{}{SIAM}.
\PrintBackRefs{\CurrentBib}

\bibitem [\protect \citeauthoryear {%
Lichtner%
\ \BBA {} Tartakovsky%
}{%
Lichtner%
\ \BBA {} Tartakovsky%
}{%
{\protect \APACyear {2003}}%
}]{%
lichtner2003stochastic}
\APACinsertmetastar {%
lichtner2003stochastic}%
\begin{APACrefauthors}%
Lichtner, P\BPBI C.%
\BCBT {}\ \BBA {} Tartakovsky, D\BPBI M.%
\end{APACrefauthors}%
\unskip\
\newblock
\APACrefYearMonthDay{2003}{}{}.
\newblock
{\BBOQ}\APACrefatitle {Stochastic analysis of effective rate constant for
  heterogeneous reactions} {Stochastic analysis of effective rate constant for
  heterogeneous reactions}.{\BBCQ}
\newblock
\APACjournalVolNumPages{Stochastic Environmental Research and Risk
  Assessment}{17}{6}{419--429}.
\PrintBackRefs{\CurrentBib}

\bibitem [\protect \citeauthoryear {%
Lucia%
, Beran%
\BCBL {}\ \BBA {} Silva%
}{%
Lucia%
\ \protect \BOthers {.}}{%
{\protect \APACyear {2004}}%
}]{%
lucia2004reduced}
\APACinsertmetastar {%
lucia2004reduced}%
\begin{APACrefauthors}%
Lucia, D\BPBI J.%
, Beran, P\BPBI S.%
\BCBL {}\ \BBA {} Silva, W\BPBI A.%
\end{APACrefauthors}%
\unskip\
\newblock
\APACrefYearMonthDay{2004}{}{}.
\newblock
{\BBOQ}\APACrefatitle {Reduced-order modeling: new approaches for computational
  physics} {Reduced-order modeling: new approaches for computational
  physics}.{\BBCQ}
\newblock
\APACjournalVolNumPages{Progress in Aerospace Sciences}{40}{1-2}{51--117}.
\PrintBackRefs{\CurrentBib}

\bibitem [\protect \citeauthoryear {%
Maday%
\ \BBA {} Mula%
}{%
Maday%
\ \BBA {} Mula%
}{%
{\protect \APACyear {2013}}%
}]{%
maday2013generalized}
\APACinsertmetastar {%
maday2013generalized}%
\begin{APACrefauthors}%
Maday, Y.%
\BCBT {}\ \BBA {} Mula, O.%
\end{APACrefauthors}%
\unskip\
\newblock
\APACrefYearMonthDay{2013}{}{}.
\newblock
{\BBOQ}\APACrefatitle {A generalized empirical interpolation method:
  application of reduced basis techniques to data assimilation} {A generalized
  empirical interpolation method: application of reduced basis techniques to
  data assimilation}.{\BBCQ}
\newblock
\BIn{} \APACrefbtitle {Analysis and Numerics of Partial Differential Equations}
  {Analysis and numerics of partial differential equations}\ (\BPGS\ 221--235).
\newblock
\APACaddressPublisher{}{Springer}.
\PrintBackRefs{\CurrentBib}

\bibitem [\protect \citeauthoryear {%
Naghibi%
, Pourghasemi%
\BCBL {}\ \BBA {} Dixon%
}{%
Naghibi%
\ \protect \BOthers {.}}{%
{\protect \APACyear {2016}}%
}]{%
naghibi2016gis}
\APACinsertmetastar {%
naghibi2016gis}%
\begin{APACrefauthors}%
Naghibi, S\BPBI A.%
, Pourghasemi, H\BPBI R.%
\BCBL {}\ \BBA {} Dixon, B.%
\end{APACrefauthors}%
\unskip\
\newblock
\APACrefYearMonthDay{2016}{}{}.
\newblock
{\BBOQ}\APACrefatitle {{GIS}-based groundwater potential mapping using boosted
  regression tree, classification and regression tree, and random forest
  machine learning models in {Iran}} {{GIS}-based groundwater potential mapping
  using boosted regression tree, classification and regression tree, and random
  forest machine learning models in {Iran}}.{\BBCQ}
\newblock
\APACjournalVolNumPages{Environmental Monitoring and Assessment}{188}{1}{44}.
\PrintBackRefs{\CurrentBib}

\bibitem [\protect \citeauthoryear {%
Neuman%
\ \BBA {} Tartakovsky%
}{%
Neuman%
\ \BBA {} Tartakovsky%
}{%
{\protect \APACyear {2009}}%
}]{%
neuman2009perspective}
\APACinsertmetastar {%
neuman2009perspective}%
\begin{APACrefauthors}%
Neuman, S\BPBI P.%
\BCBT {}\ \BBA {} Tartakovsky, D\BPBI M.%
\end{APACrefauthors}%
\unskip\
\newblock
\APACrefYearMonthDay{2009}{}{}.
\newblock
{\BBOQ}\APACrefatitle {Perspective on theories of {non-Fickian} transport in
  heterogeneous media} {Perspective on theories of {non-Fickian} transport in
  heterogeneous media}.{\BBCQ}
\newblock
\APACjournalVolNumPages{Advances in Water Resources}{32}{5}{670--680}.
\PrintBackRefs{\CurrentBib}

\bibitem [\protect \citeauthoryear {%
Opitz%
\ \BBA {} Maclin%
}{%
Opitz%
\ \BBA {} Maclin%
}{%
{\protect \APACyear {1999}}%
}]{%
opitz1999popular}
\APACinsertmetastar {%
opitz1999popular}%
\begin{APACrefauthors}%
Opitz, D.%
\BCBT {}\ \BBA {} Maclin, R.%
\end{APACrefauthors}%
\unskip\
\newblock
\APACrefYearMonthDay{1999}{}{}.
\newblock
{\BBOQ}\APACrefatitle {Popular ensemble methods: {An} empirical study} {Popular
  ensemble methods: {An} empirical study}.{\BBCQ}
\newblock
\APACjournalVolNumPages{Journal of Artificial Intelligence
  Research}{11}{}{169--198}.
\PrintBackRefs{\CurrentBib}

\bibitem [\protect \citeauthoryear {%
Pau%
, Zhang%
\BCBL {}\ \BBA {} Finsterle%
}{%
Pau%
\ \protect \BOthers {.}}{%
{\protect \APACyear {2013}}%
}]{%
pau2013reduced}
\APACinsertmetastar {%
pau2013reduced}%
\begin{APACrefauthors}%
Pau, G\BPBI S\BPBI H.%
, Zhang, Y.%
\BCBL {}\ \BBA {} Finsterle, S.%
\end{APACrefauthors}%
\unskip\
\newblock
\APACrefYearMonthDay{2013}{}{}.
\newblock
{\BBOQ}\APACrefatitle {Reduced order models for many-query subsurface flow
  applications} {Reduced order models for many-query subsurface flow
  applications}.{\BBCQ}
\newblock
\APACjournalVolNumPages{Computational Geosciences}{17}{4}{705--721}.
\PrintBackRefs{\CurrentBib}

\bibitem [\protect \citeauthoryear {%
Pedregosa%
\ \protect \BOthers {.}}{%
Pedregosa%
\ \protect \BOthers {.}}{%
{\protect \APACyear {2011}}%
}]{%
scikit-learn}
\APACinsertmetastar {%
scikit-learn}%
\begin{APACrefauthors}%
Pedregosa, F.%
, Varoquaux, G.%
, Gramfort, A.%
, Michel, V.%
, Thirion, B.%
, Grisel, O.%
\BDBL {}Duchesnay, E.%
\end{APACrefauthors}%
\unskip\
\newblock
\APACrefYearMonthDay{2011}{}{}.
\newblock
{\BBOQ}\APACrefatitle {Scikit-learn: Machine Learning in {P}ython}
  {Scikit-learn: Machine learning in {P}ython}.{\BBCQ}
\newblock
\APACjournalVolNumPages{Journal of Machine Learning
  Research}{12}{}{2825--2830}.
\PrintBackRefs{\CurrentBib}

\bibitem [\protect \citeauthoryear {%
Pinkus%
}{%
Pinkus%
}{%
{\protect \APACyear {1999}}%
}]{%
pinkus1999approximation}
\APACinsertmetastar {%
pinkus1999approximation}%
\begin{APACrefauthors}%
Pinkus, A.%
\end{APACrefauthors}%
\unskip\
\newblock
\APACrefYearMonthDay{1999}{}{}.
\newblock
{\BBOQ}\APACrefatitle {Approximation theory of the {MLP} model in neural
  networks} {Approximation theory of the {MLP} model in neural
  networks}.{\BBCQ}
\newblock
\APACjournalVolNumPages{Acta Numerica}{8}{1}{143--195}.
\PrintBackRefs{\CurrentBib}

\bibitem [\protect \citeauthoryear {%
Polikar%
}{%
Polikar%
}{%
{\protect \APACyear {2006}}%
}]{%
polikar2006ensemble}
\APACinsertmetastar {%
polikar2006ensemble}%
\begin{APACrefauthors}%
Polikar, R.%
\end{APACrefauthors}%
\unskip\
\newblock
\APACrefYearMonthDay{2006}{}{}.
\newblock
{\BBOQ}\APACrefatitle {Ensemble based systems in decision making} {Ensemble
  based systems in decision making}.{\BBCQ}
\newblock
\APACjournalVolNumPages{IEEE Circuits and Systems Magazine}{6}{3}{21--45}.
\PrintBackRefs{\CurrentBib}

\bibitem [\protect \citeauthoryear {%
Qin%
, Chen%
, Jakeman%
\BCBL {}\ \BBA {} Xiu%
}{%
Qin%
\ \protect \BOthers {.}}{%
{\protect \APACyear {2020}}%
}]{%
qin2020data}
\APACinsertmetastar {%
qin2020data}%
\begin{APACrefauthors}%
Qin, T.%
, Chen, Z.%
, Jakeman, J.%
\BCBL {}\ \BBA {} Xiu, D.%
\end{APACrefauthors}%
\unskip\
\newblock
\APACrefYearMonthDay{2020}{}{}.
\newblock
{\BBOQ}\APACrefatitle {Data-driven learning of non-autonomous systems}
  {Data-driven learning of non-autonomous systems}.{\BBCQ}
\newblock
\APACjournalVolNumPages{arXiv preprint arXiv:2006.02392}{}{}{}.
\PrintBackRefs{\CurrentBib}

\bibitem [\protect \citeauthoryear {%
Qin%
, Wu%
\BCBL {}\ \BBA {} Xiu%
}{%
Qin%
\ \protect \BOthers {.}}{%
{\protect \APACyear {2019}}%
}]{%
qin2019data}
\APACinsertmetastar {%
qin2019data}%
\begin{APACrefauthors}%
Qin, T.%
, Wu, K.%
\BCBL {}\ \BBA {} Xiu, D.%
\end{APACrefauthors}%
\unskip\
\newblock
\APACrefYearMonthDay{2019}{}{}.
\newblock
{\BBOQ}\APACrefatitle {Data driven governing equations approximation using deep
  neural networks} {Data driven governing equations approximation using deep
  neural networks}.{\BBCQ}
\newblock
\APACjournalVolNumPages{Journal of Computational Physics}{395}{}{620--635}.
\PrintBackRefs{\CurrentBib}

\bibitem [\protect \citeauthoryear {%
Rasmussen%
}{%
Rasmussen%
}{%
{\protect \APACyear {2003}}%
}]{%
rasmussen2003gaussian}
\APACinsertmetastar {%
rasmussen2003gaussian}%
\begin{APACrefauthors}%
Rasmussen, C\BPBI E.%
\end{APACrefauthors}%
\unskip\
\newblock
\APACrefYearMonthDay{2003}{}{}.
\newblock
{\BBOQ}\APACrefatitle {Gaussian processes in machine learning} {Gaussian
  processes in machine learning}.{\BBCQ}
\newblock
\BIn{} \APACrefbtitle {Summer School on Machine Learning} {Summer school on
  machine learning}\ (\BPGS\ 63--71).
\PrintBackRefs{\CurrentBib}

\bibitem [\protect \citeauthoryear {%
Razavi%
, Tolson%
\BCBL {}\ \BBA {} Burn%
}{%
Razavi%
\ \protect \BOthers {.}}{%
{\protect \APACyear {2012}}%
}]{%
razavi2012review}
\APACinsertmetastar {%
razavi2012review}%
\begin{APACrefauthors}%
Razavi, S.%
, Tolson, B\BPBI A.%
\BCBL {}\ \BBA {} Burn, D\BPBI H.%
\end{APACrefauthors}%
\unskip\
\newblock
\APACrefYearMonthDay{2012}{}{}.
\newblock
{\BBOQ}\APACrefatitle {Review of surrogate modeling in water resources} {Review
  of surrogate modeling in water resources}.{\BBCQ}
\newblock
\APACjournalVolNumPages{Water Resources Research}{48}{7}{}.
\PrintBackRefs{\CurrentBib}

\bibitem [\protect \citeauthoryear {%
Rokach%
}{%
Rokach%
}{%
{\protect \APACyear {2010}}%
}]{%
rokach2010ensemble}
\APACinsertmetastar {%
rokach2010ensemble}%
\begin{APACrefauthors}%
Rokach, L.%
\end{APACrefauthors}%
\unskip\
\newblock
\APACrefYearMonthDay{2010}{}{}.
\newblock
{\BBOQ}\APACrefatitle {Ensemble-based classifiers} {Ensemble-based
  classifiers}.{\BBCQ}
\newblock
\APACjournalVolNumPages{Artificial Intelligence Review}{33}{1-2}{1--39}.
\PrintBackRefs{\CurrentBib}

\bibitem [\protect \citeauthoryear {%
Rowley%
}{%
Rowley%
}{%
{\protect \APACyear {2005}}%
}]{%
rowley2005model}
\APACinsertmetastar {%
rowley2005model}%
\begin{APACrefauthors}%
Rowley, C\BPBI W.%
\end{APACrefauthors}%
\unskip\
\newblock
\APACrefYearMonthDay{2005}{}{}.
\newblock
{\BBOQ}\APACrefatitle {Model reduction for fluids, using balanced proper
  orthogonal decomposition} {Model reduction for fluids, using balanced proper
  orthogonal decomposition}.{\BBCQ}
\newblock
\APACjournalVolNumPages{International Journal of Bifurcation and
  Chaos}{15}{03}{997--1013}.
\PrintBackRefs{\CurrentBib}

\bibitem [\protect \citeauthoryear {%
Rutqvist%
, Zheng%
, Chen%
, Liu%
\BCBL {}\ \BBA {} Birkholzer%
}{%
Rutqvist%
\ \protect \BOthers {.}}{%
{\protect \APACyear {2014}}%
}]{%
rutqvist2014modeling}
\APACinsertmetastar {%
rutqvist2014modeling}%
\begin{APACrefauthors}%
Rutqvist, J.%
, Zheng, L.%
, Chen, F.%
, Liu, H\BHBI H.%
\BCBL {}\ \BBA {} Birkholzer, J.%
\end{APACrefauthors}%
\unskip\
\newblock
\APACrefYearMonthDay{2014}{}{}.
\newblock
{\BBOQ}\APACrefatitle {Modeling of coupled thermo-hydro-mechanical processes
  with links to geochemistry associated with bentonite-backfilled repository
  tunnels in clay formations} {Modeling of coupled thermo-hydro-mechanical
  processes with links to geochemistry associated with bentonite-backfilled
  repository tunnels in clay formations}.{\BBCQ}
\newblock
\APACjournalVolNumPages{Rock Mechanics and Rock Engineering}{47}{1}{167--186}.
\PrintBackRefs{\CurrentBib}

\bibitem [\protect \citeauthoryear {%
Saridakis%
\ \BBA {} Dentsoras%
}{%
Saridakis%
\ \BBA {} Dentsoras%
}{%
{\protect \APACyear {2008}}%
}]{%
saridakis2008soft}
\APACinsertmetastar {%
saridakis2008soft}%
\begin{APACrefauthors}%
Saridakis, K\BPBI M.%
\BCBT {}\ \BBA {} Dentsoras, A\BPBI J.%
\end{APACrefauthors}%
\unskip\
\newblock
\APACrefYearMonthDay{2008}{}{}.
\newblock
{\BBOQ}\APACrefatitle {Soft computing in engineering design--{A} review} {Soft
  computing in engineering design--{A} review}.{\BBCQ}
\newblock
\APACjournalVolNumPages{Advanced Engineering Informatics}{22}{2}{202--221}.
\PrintBackRefs{\CurrentBib}

\bibitem [\protect \citeauthoryear {%
Schmid%
}{%
Schmid%
}{%
{\protect \APACyear {2010}}%
}]{%
schmid2010dynamic}
\APACinsertmetastar {%
schmid2010dynamic}%
\begin{APACrefauthors}%
Schmid, P\BPBI J.%
\end{APACrefauthors}%
\unskip\
\newblock
\APACrefYearMonthDay{2010}{}{}.
\newblock
{\BBOQ}\APACrefatitle {Dynamic mode decomposition of numerical and experimental
  data} {Dynamic mode decomposition of numerical and experimental data}.{\BBCQ}
\newblock
\APACjournalVolNumPages{Journal of Fluid Mechanics}{656}{}{5--28}.
\PrintBackRefs{\CurrentBib}

\bibitem [\protect \citeauthoryear {%
Schmit~Jr%
\ \BBA {} Farshi%
}{%
Schmit~Jr%
\ \BBA {} Farshi%
}{%
{\protect \APACyear {1974}}%
}]{%
schmit1974some}
\APACinsertmetastar {%
schmit1974some}%
\begin{APACrefauthors}%
Schmit~Jr, L.%
\BCBT {}\ \BBA {} Farshi, B.%
\end{APACrefauthors}%
\unskip\
\newblock
\APACrefYearMonthDay{1974}{}{}.
\newblock
{\BBOQ}\APACrefatitle {Some approximation concepts for structural synthesis}
  {Some approximation concepts for structural synthesis}.{\BBCQ}
\newblock
\APACjournalVolNumPages{AIAA Journal}{12}{5}{692--699}.
\PrintBackRefs{\CurrentBib}

\bibitem [\protect \citeauthoryear {%
Simpson%
, Poplinski%
, Koch%
\BCBL {}\ \BBA {} Allen%
}{%
Simpson%
\ \protect \BOthers {.}}{%
{\protect \APACyear {2001}}%
}]{%
simpson2001metamodels}
\APACinsertmetastar {%
simpson2001metamodels}%
\begin{APACrefauthors}%
Simpson, T\BPBI W.%
, Poplinski, J.%
, Koch, P\BPBI N.%
\BCBL {}\ \BBA {} Allen, J\BPBI K.%
\end{APACrefauthors}%
\unskip\
\newblock
\APACrefYearMonthDay{2001}{}{}.
\newblock
{\BBOQ}\APACrefatitle {Metamodels for computer-based engineering design: survey
  and recommendations} {Metamodels for computer-based engineering design:
  survey and recommendations}.{\BBCQ}
\newblock
\APACjournalVolNumPages{Engineering with Computers}{17}{2}{129--150}.
\PrintBackRefs{\CurrentBib}

\bibitem [\protect \citeauthoryear {%
Steefel%
\ \protect \BOthers {.}}{%
Steefel%
\ \protect \BOthers {.}}{%
{\protect \APACyear {2015}}%
}]{%
steefel2015reactive}
\APACinsertmetastar {%
steefel2015reactive}%
\begin{APACrefauthors}%
Steefel, C\BPBI I.%
, Appelo, C\BPBI A\BPBI J.%
, Arora, B.%
, Jacques, D.%
, Kalbacher, T.%
, Kolditz, O.%
\BDBL {}others%
\end{APACrefauthors}%
\unskip\
\newblock
\APACrefYearMonthDay{2015}{}{}.
\newblock
{\BBOQ}\APACrefatitle {Reactive transport codes for subsurface environmental
  simulation} {Reactive transport codes for subsurface environmental
  simulation}.{\BBCQ}
\newblock
\APACjournalVolNumPages{Computational Geosciences}{19}{3}{445--478}.
\PrintBackRefs{\CurrentBib}

\bibitem [\protect \citeauthoryear {%
Steefel%
, DePaolo%
\BCBL {}\ \BBA {} Lichtner%
}{%
Steefel%
\ \protect \BOthers {.}}{%
{\protect \APACyear {2005}}%
}]{%
steefel2005reactive}
\APACinsertmetastar {%
steefel2005reactive}%
\begin{APACrefauthors}%
Steefel, C\BPBI I.%
, DePaolo, D\BPBI J.%
\BCBL {}\ \BBA {} Lichtner, P\BPBI C.%
\end{APACrefauthors}%
\unskip\
\newblock
\APACrefYearMonthDay{2005}{}{}.
\newblock
{\BBOQ}\APACrefatitle {Reactive transport modeling: {An} essential tool and a
  new research approach for the {Earth} sciences} {Reactive transport modeling:
  {An} essential tool and a new research approach for the {Earth}
  sciences}.{\BBCQ}
\newblock
\APACjournalVolNumPages{Earth and Planetary Science
  Letters}{240}{3-4}{539--558}.
\PrintBackRefs{\CurrentBib}

\bibitem [\protect \citeauthoryear {%
Tavenard%
\ \protect \BOthers {.}}{%
Tavenard%
\ \protect \BOthers {.}}{%
{\protect \APACyear {2020}}%
}]{%
JMLR:v21:20-091}
\APACinsertmetastar {%
JMLR:v21:20-091}%
\begin{APACrefauthors}%
Tavenard, R.%
, Faouzi, J.%
, Vandewiele, G.%
, Divo, F.%
, Androz, G.%
, Holtz, C.%
\BDBL {}Woods, E.%
\end{APACrefauthors}%
\unskip\
\newblock
\APACrefYearMonthDay{2020}{}{}.
\newblock
{\BBOQ}\APACrefatitle {Tslearn, A Machine Learning Toolkit for Time Series
  Data} {Tslearn, a machine learning toolkit for time series data}.{\BBCQ}
\newblock
\APACjournalVolNumPages{Journal of Machine Learning Research}{21}{118}{1-6}.
\PrintBackRefs{\CurrentBib}

\bibitem [\protect \citeauthoryear {%
Tripathy%
\ \BBA {} Bilionis%
}{%
Tripathy%
\ \BBA {} Bilionis%
}{%
{\protect \APACyear {2018}}%
}]{%
tripathy2018deep}
\APACinsertmetastar {%
tripathy2018deep}%
\begin{APACrefauthors}%
Tripathy, R\BPBI K.%
\BCBT {}\ \BBA {} Bilionis, I.%
\end{APACrefauthors}%
\unskip\
\newblock
\APACrefYearMonthDay{2018}{}{}.
\newblock
{\BBOQ}\APACrefatitle {{Deep UQ: Learning} deep neural network surrogate models
  for high dimensional uncertainty quantification} {{Deep UQ: Learning} deep
  neural network surrogate models for high dimensional uncertainty
  quantification}.{\BBCQ}
\newblock
\APACjournalVolNumPages{Journal of Computational Physics}{375}{}{565--588}.
\PrintBackRefs{\CurrentBib}

\bibitem [\protect \citeauthoryear {%
Xiong%
\ \protect \BOthers {.}}{%
Xiong%
\ \protect \BOthers {.}}{%
{\protect \APACyear {2013}}%
}]{%
xiong2013coupled}
\APACinsertmetastar {%
xiong2013coupled}%
\begin{APACrefauthors}%
Xiong, Y.%
, Fakcharoenphol, P.%
, Winterfeld, P.%
, Zhang, R.%
, Wu, Y\BHBI S.%
\BCBL {}\ \BOthersPeriod {.}\end{APACrefauthors}%
\unskip\
\newblock
\APACrefYearMonthDay{2013}{}{}.
\newblock
{\BBOQ}\APACrefatitle {Coupled geomechanical and reactive geochemical model for
  fluid and heat flow: application for enhanced geothermal reservoir} {Coupled
  geomechanical and reactive geochemical model for fluid and heat flow:
  application for enhanced geothermal reservoir}.{\BBCQ}
\newblock
\BIn{} \APACrefbtitle {{SPE} Reservoir Characterization and Simulation
  Conference and Exhibition.} {{SPE} reservoir characterization and simulation
  conference and exhibition.}
\PrintBackRefs{\CurrentBib}

\bibitem [\protect \citeauthoryear {%
Xu%
, Sonnenthal%
, Spycher%
\BCBL {}\ \BBA {} Zheng%
}{%
Xu%
\ \protect \BOthers {.}}{%
{\protect \APACyear {2014}}%
}]{%
xu2014toughreact}
\APACinsertmetastar {%
xu2014toughreact}%
\begin{APACrefauthors}%
Xu, T.%
, Sonnenthal, E.%
, Spycher, N.%
\BCBL {}\ \BBA {} Zheng, L.%
\end{APACrefauthors}%
\unskip\
\newblock
\APACrefYearMonthDay{2014}{}{}.
\newblock
{\BBOQ}\APACrefatitle {TOUGHREACT V3. 0-OMP reference manual: A parallel
  simulation program for non-isothermal multiphase geochemical reactive
  transport} {Toughreact v3. 0-omp reference manual: A parallel simulation
  program for non-isothermal multiphase geochemical reactive transport}.{\BBCQ}
\newblock
\APACjournalVolNumPages{University of California, Berkeley}{}{}{}.
\PrintBackRefs{\CurrentBib}

\bibitem [\protect \citeauthoryear {%
Zheng%
, Rutqvist%
, Xu%
\BCBL {}\ \BBA {} Birkholzer%
}{%
Zheng%
\ \protect \BOthers {.}}{%
{\protect \APACyear {2017}}%
}]{%
zheng2017coupled}
\APACinsertmetastar {%
zheng2017coupled}%
\begin{APACrefauthors}%
Zheng, L.%
, Rutqvist, J.%
, Xu, H.%
\BCBL {}\ \BBA {} Birkholzer, J\BPBI T.%
\end{APACrefauthors}%
\unskip\
\newblock
\APACrefYearMonthDay{2017}{}{}.
\newblock
{\BBOQ}\APACrefatitle {Coupled THMC models for bentonite in an argillite
  repository for nuclear waste: Illitization and its effect on swelling stress
  under high temperature} {Coupled thmc models for bentonite in an argillite
  repository for nuclear waste: Illitization and its effect on swelling stress
  under high temperature}.{\BBCQ}
\newblock
\APACjournalVolNumPages{Engineering geology}{230}{}{118--129}.
\PrintBackRefs{\CurrentBib}

\bibitem [\protect \citeauthoryear {%
Zhou%
\ \BBA {} Tartakovsky%
}{%
Zhou%
\ \BBA {} Tartakovsky%
}{%
{\protect \APACyear {2020}}%
}]{%
zhou-2020-markov}
\APACinsertmetastar {%
zhou-2020-markov}%
\begin{APACrefauthors}%
Zhou, Z.%
\BCBT {}\ \BBA {} Tartakovsky, D\BPBI M.%
\end{APACrefauthors}%
\unskip\
\newblock
\APACrefYearMonthDay{2020}{}{}.
\newblock
{\BBOQ}\APACrefatitle {Markov chain {Monte Carlo} with neural network
  surrogates: {Application} to contaminant source identification} {Markov chain
  {Monte Carlo} with neural network surrogates: {Application} to contaminant
  source identification}.{\BBCQ}
\newblock
\APACjournalVolNumPages{Stochastic Environmental Research and Risk
  Assessment}{}{}{}.
\newblock
\begin{APACrefDOI} \doi{10.1007/s00477-020-01888-9} \end{APACrefDOI}
\PrintBackRefs{\CurrentBib}

\bibitem [\protect \citeauthoryear {%
R.~Zhu%
, Zeng%
\BCBL {}\ \BBA {} Kosorok%
}{%
R.~Zhu%
\ \protect \BOthers {.}}{%
{\protect \APACyear {2015}}%
}]{%
zhu2015reinforcement}
\APACinsertmetastar {%
zhu2015reinforcement}%
\begin{APACrefauthors}%
Zhu, R.%
, Zeng, D.%
\BCBL {}\ \BBA {} Kosorok, M\BPBI R.%
\end{APACrefauthors}%
\unskip\
\newblock
\APACrefYearMonthDay{2015}{}{}.
\newblock
{\BBOQ}\APACrefatitle {Reinforcement learning trees} {Reinforcement learning
  trees}.{\BBCQ}
\newblock
\APACjournalVolNumPages{Journal of the American Statistical
  Association}{110}{512}{1770--1784}.
\PrintBackRefs{\CurrentBib}

\bibitem [\protect \citeauthoryear {%
Y.~Zhu%
, Zabaras%
, Koutsourelakis%
\BCBL {}\ \BBA {} Perdikaris%
}{%
Y.~Zhu%
\ \protect \BOthers {.}}{%
{\protect \APACyear {2019}}%
}]{%
zhu2019physics}
\APACinsertmetastar {%
zhu2019physics}%
\begin{APACrefauthors}%
Zhu, Y.%
, Zabaras, N.%
, Koutsourelakis, P\BHBI S.%
\BCBL {}\ \BBA {} Perdikaris, P.%
\end{APACrefauthors}%
\unskip\
\newblock
\APACrefYearMonthDay{2019}{}{}.
\newblock
{\BBOQ}\APACrefatitle {Physics-constrained deep learning for high-dimensional
  surrogate modeling and uncertainty quantification without labeled data}
  {Physics-constrained deep learning for high-dimensional surrogate modeling
  and uncertainty quantification without labeled data}.{\BBCQ}
\newblock
\APACjournalVolNumPages{Journal of Computational Physics}{394}{}{56--81}.
\PrintBackRefs{\CurrentBib}

\end{thebibliography}

\end{document}